\documentclass[11pt,reqno]{article}

\usepackage{amsfonts}
\usepackage{amsmath}
\usepackage{amsthm}
\usepackage{amssymb}
\usepackage{mathrsfs}
\usepackage{amstext}

\usepackage[dvips]{graphicx}
\usepackage{subfigure}
\usepackage{algorithm}
\usepackage{algorithmic}

\font\msbm=msbm10

\numberwithin{equation}{section}

\theoremstyle{plain}
\newtheorem{Theorem}{Theorem}[section]
\newtheorem{lemma}[Theorem]{Lemma}

\newtheorem{example}[Theorem]{Example}

\newtheorem{remark}[Theorem]{Remark}
\def\mathbb#1{\hbox{\msbm{#1}}}
\newcommand{\N}{{\mathbb{N}}}
\newcommand{\R}{{\mathbb{R}}}
\newcommand{\Z}{{\mathbb{Z}}}
\newcommand{\C}{{\mathbb{C}}}

\renewcommand{\P}{{\mathbb{P}}}
\newcommand{\E}{{\mathbb{E}}}

\newcommand{\beq}{\begin{eqnarray}}
\newcommand{\eeq}{\end{eqnarray}}
\newcommand{\beqn}{\begin{eqnarray*}}
\newcommand{\eeqn}{\end{eqnarray*}}

\newcommand{\ff}{{\mathrm f}}
\newcommand{\F}{{\cal F}}

\newcommand{\supp}{\operatorname{supp}}

\newcommand{\spn}{\operatorname{span}}
\newcommand{\re}{\operatorname{Re}}
\newcommand{\im}{\operatorname{Im}}

\renewcommand{\qed}{\rule{2.5mm}{2.5mm}}
\newenvironment{Proof}{\noindent {\bf\underline{Proof:}}}{\hspace*{\fill}\qed\vskip1em}

\begin{document}
\title{Random Sampling of Sparse Trigonometric Polynomials II - Orthogonal Matching Pursuit versus Basis Pursuit}

\author{Stefan Kunis, Holger Rauhut}
\date{April 19, 2006; revised \today}  
%NuHAG, Faculty of Mathematics, University of Vienna\\
%Nordbergstrasse 15, A-1090 Wien, Austria\\
%rauhut@ma.tum.de

\maketitle

\begin{abstract}
We investigate the problem of reconstructing sparse multivariate 
trigonometric polynomials from few randomly taken samples by Basis Pursuit and 
greedy algorithms such as Orthogonal Matching Pursuit (OMP) and Thresholding. 
While recovery by Basis Pursuit has recently been studied by several authors,
we provide theoretical results 
on the success probability of reconstruction via Thresholding and OMP for
both a continuous and a discrete 
probability model for the sampling points.
We present numerical experiments, which indicate that
usually Basis Pursuit is significantly slower than
greedy algorithms, while the recovery rates are very similar. 
\end{abstract}
\vspace{0.5cm}


%94A20: Sampling Theory
%42A05: Trigonometric Polynomials%15A52: Random matrices
%05A18: Enumerative Combinatorics: Partitions of Sets
%90C05: linear programming
%90C25: convex programming

%-------------------------------------------------------------------------------------------------------
\section{Introduction}

In the last two years the rapidly growing field of
compressed sensing has attracted much attention 
\cite{badadewa06,CRT1,CRT2,CT,CT2,Donoho1,Tropp,Rau,RV06}.
Its basic idea is that sparse or compressible
signals can be reconstructed from vastly incomplete
non-adaptive information.
By ``sparse'' we mean that a vector has only few non-zero coefficients,
while ``compressible'' expresses that a vector can be well-approximated
by a sparse one.

Previous work on this topic includes the reconstruction of Fourier
coefficients from samples taken randomly on a lattice 
by the {\em Basis Pursuit} (BP) principle \cite{CRT1,CT}.
This consists in minimizing the $\ell^1$-norm of the 
Fourier coefficients subject to the condition that the corresponding
trigonometric polynomial matches the sampling points.
Indeed, it was proven by Cand{\`e}s, Romberg and Tao in \cite{CRT1} in the setting
of the discrete Fourier transform that this scheme recovers the coefficients
exactly with high probability provided the number of samples is high enough
compared to the sparsity, i.e., the number of non-zero coefficients.
This result has been generalized by the second author of the present paper in
\cite{Rau} for the case of samples taken uniformly at random from the cube
$[0,2\pi]^d$.

Another line of research
suggests greedy methods such as (Orthogonal) Matching Pursuit (OMP) and 
Thresholding for sparse reconstruction tasks 
\cite{MaZh,Tropp, Tropp_Greed, Tropp_Thesis}.
OMP and Thresholding are conceptually simple to implement and
potentially faster than BP.
In particular, they may easily take into account fast algorithms for 
multiplication with the involved matrices, while most standard software for convex
optimization \cite{BV,MOSEK} does not allow this.

This paper is devoted to the theoretical and numerical investigation and
comparison of Thresholding, OMP and BP for the recovery of sparse trigonometric
polynomials from randomly taken samples.
Our theoretical results indicate that indeed all the
methods are suitable for this task. The novelty in the present paper
is a performance analysis for OMP and Thresholding, 
although the theoretical achievements for OMP are only partial so far. 
In contrast to BP, the greedy algorithms give only a non-uniform guarantee of 
recovery at a sufficiently small ratio of the number of samples to the 
sparsity. This means that they guarantee recovery with high probability 
only for the given trigonometric polynomial, while
BP can actually guarantee recovery of {\em all} sufficiently sparse
trigonometric polynomials from a single sampling set.

In practice however, a non-uniform guarantee might be sufficient. Indeed,
our numerical experiments suggest that OMP even slightly outperforms BP 
on generic signals with respect to reconstruction rate. Considering
that greedy algorithms are usually 
significantly faster than BP one would probably use
OMP for most applications despite its lack of uniformity.

For related work on this topic, also known as compressed sensing,
we refer to 
\cite{badadewa06,CRT1,CT,CRT2,CR2,CR,Donoho1,DTsaig,Donoho2,DT,Tropp,RV06} 
and the references therein.
For more information on sampling of (not necessarily sparse) trigonometric
polynomials in a probabilistic setting the reader may consult \cite{BG,BP,grrapo06}.

\begin{figure}[h]
  \centering
  \subfigure[Sparse coefficient vector.]   
  {\includegraphics[width=0.45\textwidth]{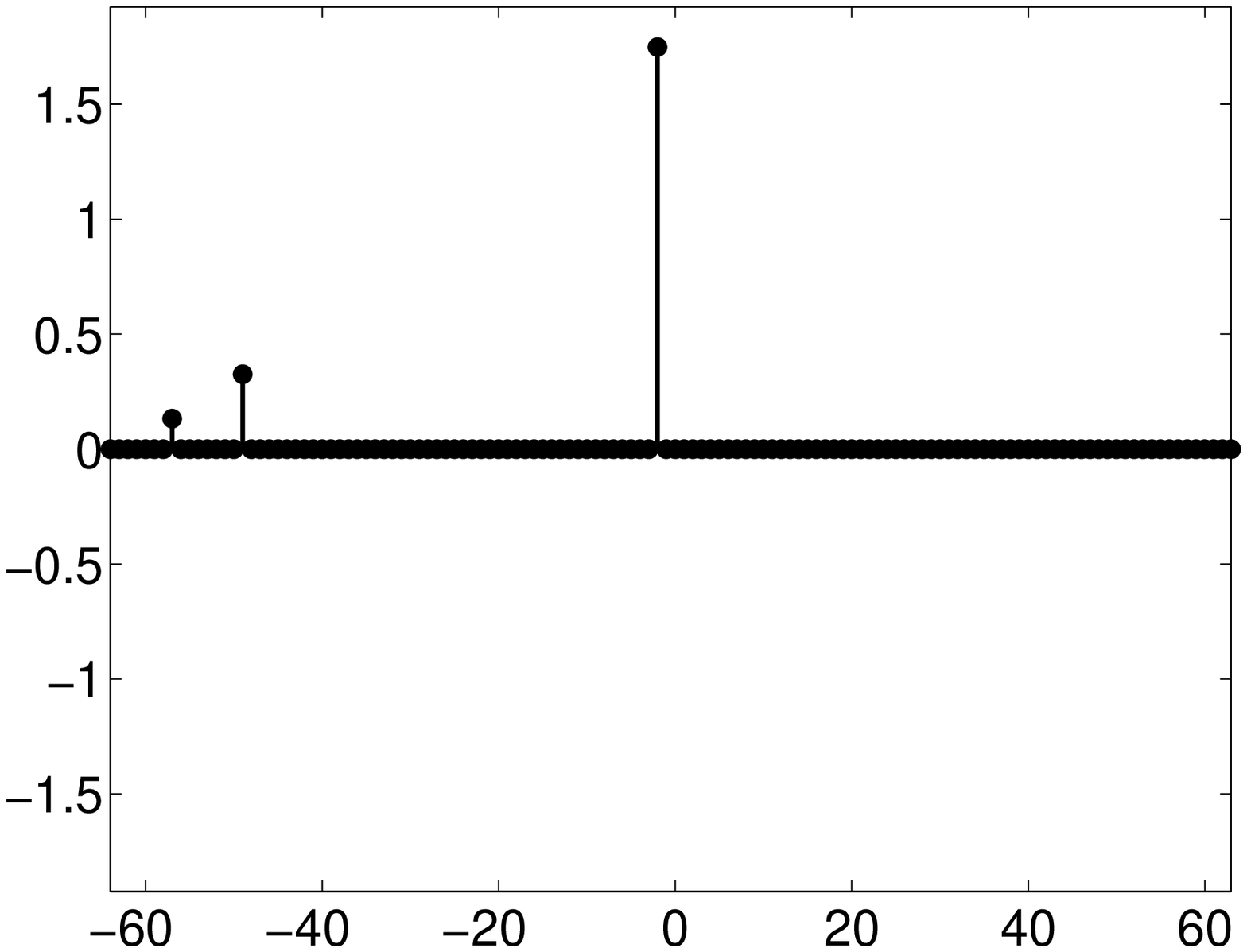}}\hfill
  \subfigure[Trigonometric polynomial and a few samples.]
  {\includegraphics[width=0.45\textwidth]{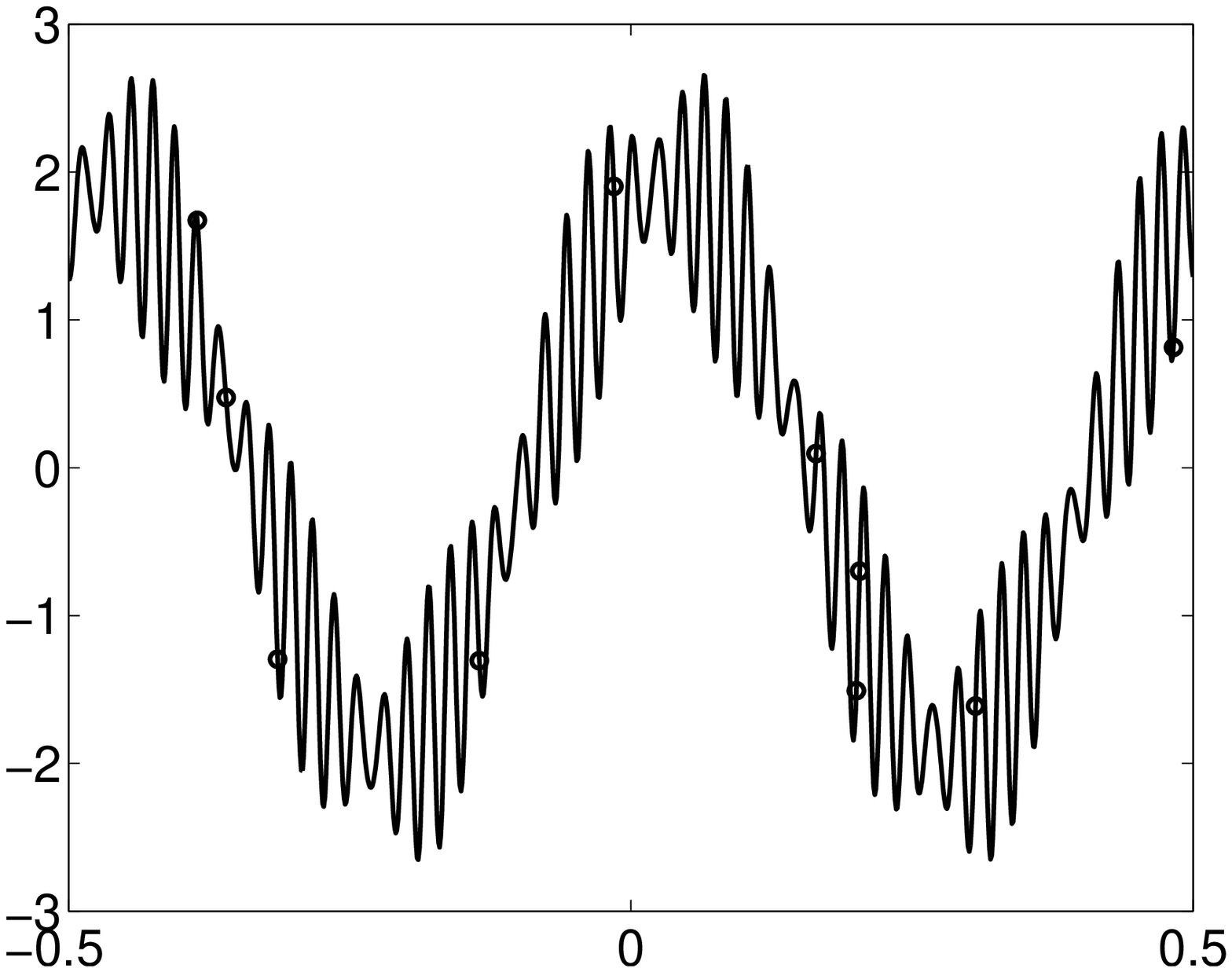}}\\
  \caption{Sparse vector of Fourier coefficients and the corresponding
    trigonometric polynomial (real part).
    After sampling at a few randomly chosen points, Orthogonal Matching
    Pursuit (OMP), i.e., Algorithm \ref{algo:omp}, as well as the
    Basis Pursuit principle recover the coefficient vector perfectly with high
    probability.\label{fig:1}}
\end{figure}

The paper is organized as follows: After introducing the necessary notation,
including the Orthogonal Matching Pursuit algorithm, we first review
known results for Basis Pursuit. Then we present our main
result concerning Thresholding and OMP.  
Based on the coherence parameter we also provide uniform reconstruction
results for Thresholding and OMP, cf.~Subsection \ref{coherence}.

In Section \ref{sect:proofs} all proofs of the obtained results are given.
Section \ref{sect:num} presents extensive 
numerical experiments.
Finally, Section 5 makes conclusions and discusses possible future work.

\section{Main Results}
\subsection{The Setting}

For some finite subset $\Gamma \subset \Z^d$, $d \in \N$, we let $\Pi_\Gamma$
denote the space of all trigonometric polynomials in dimension $d$ whose
coefficients are supported on $\Gamma$.
Clearly, an element $f$ of $\Pi_\Gamma$ is of the form 
\[
f(x) \,=\, \sum_{k \in \Gamma} c_k e^{i k\cdot x}, \qquad 
x \in [0,2\pi]^d, 
\]
with some Fourier coefficients $c_k \in \C$.
The dimension of $\Pi_\Gamma$ will be denoted by $D:= |\Gamma|$.
Taking $\Gamma = \{-q,-q+1,\hdots,q\}^d$ yields the space $\Pi_\Gamma =
\Pi_q^d$  of all trigonometric polynomials of maximal order $q$. 

We will mainly deal with ``sparse'' trigonometric polynomials, i.e., we assume
that the sequence of coefficients $c_k$ is supported only on a set $T$, which
is much smaller than $\Gamma$.
However, a priori nothing is known about $T$ except for its maximum size.
Thus, it is useful to introduce the set $\Pi_\Gamma(M) \subset \Pi_\Gamma$ of
all trigonometric polynomials whose Fourier coefficients are supported on a
set  $T \subset \Gamma$ satisfying $|T| \leq M$.
Note that 
\[
\Pi_\Gamma(M)\,=\, \bigcup_{T\subset \Gamma, |T|\leq M} \Pi_T
\]
is not a linear space.

Our aim is to sample a trigonometric polynomial $f$ of $\Pi_\Gamma(M)$ at $N$
points $x_1,\hdots, x_N\in [0,2\pi]^d$ and try to reconstruct $f$ from these
samples.
If for some $m\in\N$ the sampling points are located on the grid 
\[
\frac{2\pi}{m} \Z_m^d \,=\, \left\{0,\frac{2\pi}{m},\hdots, \frac{2\pi(m-1)}{m}
\right\}
\] 
then this problem 
can also be interpreted as reconstructing a sparse vector from partial
information on its discrete Fourier transform.

Basis Pursuit consists in solving the following $\ell^1$-minimization problem
\begin{equation}\label{BP}
\min \|d\|_1:=\sum_{k \in \Gamma} |d_k| \quad \mbox{ subject to } \quad \sum_{k \in \Gamma } 
d_k e^{2\pi i k\cdot x_j} = f\left(x_j\right),\quad j=1,\hdots,N.
\end{equation}
This task can be performed with convex optimization techniques \cite{BV}.
For real-valued coefficients (\ref{BP}) can be reformulated as
a linear program while for complex-valued coefficients we obtain
a second order cone program. For both kind of problems standard software
exists, such as MOSEK \cite{MOSEK} or CVX \cite{GrBo_CVX} 
(internally using SeDuMi \cite{SeDuMi}) and
since recently also L1MAGIC \cite{roca_l1MAGIC} 
(only for real-valued coefficients).

As an alternative to BP we also use greedy algorithms 
to recover the Fourier  coefficients of $f$ from few samples.
In particular, we study Orthogonal Matching Pursuit (Algorithm \ref{algo:omp})
as well as
the very simple Thresholding algorithm (Algorithm \ref{algo:thresh}).
We need to introduce some notation.
Let $X=\{x_1,\hdots,x_N\}$ be the set of (random) sampling points.  
We denote by $\F_X$ the $N \times D$ matrix (recall that $D = |\Gamma|$) with
entries
\begin{equation}\label{def_FX}
  (\F_X)_{j,k} \,=\, e^{ik \cdot x_j}, \quad 1\leq j \leq N,\, k \in \Gamma.
\end{equation}
Then clearly, $f(x_j) = (\F_X c)_j$ if $c$ is the vector of Fourier
coefficients of $f$. Let $\phi_k$ denote the $k$-th column of $\F_X$, i.e.,
\[
\phi_k \,=\, \left(\begin{matrix} e^{ik\cdot x_1} \\ \vdots \\ e^{ik\cdot x_N}
  \end{matrix}\right),
\]
so $\F_X = (\phi_{k_1}|\phi_{k_2}|\hdots|\phi_{k_D})$.
By 
\begin{equation}\label{def_FTX}
  (\F_{TX})_{j,k} \,=\, e^{ik \cdot x_j}, \quad 1\leq j \leq N,\, k \in T.
\end{equation}
we denote the restriction of $\F_X$ to sequences supported only on $T$.
Furthermore, let $\langle \cdot, \cdot \rangle$ be the usual Euclidean
scalar product and $\|\cdot\|_2$ the associated norm.
We have $\| \phi_k \|_2 = \sqrt{N}$ for all $k \in \Gamma$, i.e., all the
columns of $\F_X$ have the same $\ell^2$-norm.
We postpone a detailed discussion on the implementation of Algorithm
\ref{algo:omp} to Section \ref{sect:num}.

\begin{algorithm}[t]
  \caption{OMP\label{algo:omp}}
  \begin{tabular}{ll}
    Input:    & sampling set $X\subset [0,2\pi]^d$, sampling vector
                $\ff:=(f(x_j))_{j=1}^N$, set $\Gamma\subset\Z^d$.\\
    Optional: & maximum allowed sparsity $M$ or residual tolerance
                $\varepsilon$.\\[1ex]
  \end{tabular}
  \begin{algorithmic}[1]
    \STATE Set $s=0$, the residual vector $r_0 = \ff$, and the index set
           $T_0 = \emptyset$.
    \REPEAT
    \STATE Set $s=s+1$.
    \STATE Find $k_s=\arg\max_{k\in\Gamma} |\langle r_{s-1},\phi_k\rangle|$
    and augment $T_s = T_{s-1} \cup \{k_s\}$.
    \STATE Project onto $\spn\{ \phi_k, k\in T_s\}$ by solving the least squares problem
    \begin{equation*}
      \left\|\F_{T_s X} d_s - \ff\right\|_2\stackrel{d_s}{\rightarrow}\min.
    \end{equation*}
    \STATE Compute the new residual $r_s =  \ff - \F_{T_s X} d_s$.
    \UNTIL{$s=M$ or $\|r_s\|\le\varepsilon$}
    \STATE Set $T = T_s$, the non-zeros of the vector $c$ are given by
    $(c_k)_{k\in T}=d_s$.
  \end{algorithmic}
  \begin{tabular}{ll}
    \\[1ex]
    Output:   & vector of coefficients $c$ %$(c_k)_{k\in \Gamma}$ 
                and its support $T$.
  \end{tabular}
\end{algorithm}

\begin{algorithm}[t]
\caption{Thresholding\label{algo:thresh}}
\begin{tabular}{ll}
Input:    & sampling set $X\subset [0,2\pi]^d$, sampling vector
                $\ff:=(f(x_j))_{j=1}^N$, set $\Gamma\subset\Z^d$,\\
      & maximum allowed sparsity $M$
\end{tabular}\\[2ex]
\begin{algorithmic}[1]
\STATE Find the indices $T \subset \Gamma$ corresponding to the $M$ largest
inner products $\{|\langle \ff, \phi_k \rangle |\}_{k \in \Gamma}$.
\STATE Project onto $\spn\{ \phi_k, k\in T_s\}$ by solving the least 
squares problem
    \begin{equation*}
      \left\|\F_{T X} d - \ff\right\|_2\stackrel{d}{\rightarrow}\min.
    \end{equation*}
\STATE The non-zero entries of the vector $c$ are given by $(c_k)_{k \in T} = d$.
\end{algorithmic}
\begin{tabular}{ll}
\\[2ex]
Output: & vector of coefficients $c$ and its support $T$
\end{tabular}
\end{algorithm}

Of course, the hope is that running OMP or BP on samples of some $f \in
\Pi_\Gamma(M)$ will recover its Fourier coefficients.
To analyze the performance of the algorithms we will use 
two probabilistic models
for the sampling points (one for the continuous case and the other one for
the discrete Fourier transform case). This random modeling 
can be understood 
in the sense that the sampling set $X$ is 'generic': reconstruction is allowed
to fail for certain choices of $X$, as long as the probability of encountering
such a pathological case is very small. 

We will work with the following two probability models for the sampling points:
\begin{itemize}
\item[(1)] We assume that the sampling points $x_1,\hdots,x_N$ 
are independent random variables having the uniform distribution on $[0,2\pi]^d$.
Obviously, the cardinality of the sampling set $X=\{ x_1,\hdots,x_N \}$ 
equals the number of samples
$N$ with probability $1$.
\item[(2)] We suppose that the sampling points
$x_1,\hdots,x_N$ have the uniform distribution on the finite set 
$\frac{2\pi}{m} \Z_m^d$ for some $m \in \N \setminus \{1\}$. 
It will then always be assumed 
implicitly that $\Gamma \subset \Z_m^d$. % when we discuss this second model. 
Clearly, our second model aims at studying the problem of 
reconstructing sparse vectors from (partial) information on
its discrete Fourier transform.

Observe that it happens with non-zero probability that some of the sampling
points coincide, so the cardinality of the sampling set $X =
\{x_1,\hdots,x_N\}$ might be smaller than $N$.
However, for small $N \ll D$ this effect will occur rather rarely, and
also does not harm our theoretical analysis.
\end{itemize}

We will often refer to the first model as the ``continuous model'' while
the second will be called the ``discrete model''.
It turns out that one can treat both probability models in parallel. 

\subsection{Previous results for Basis Pursuit}
\label{Sec:BP}

Based on ideas due to Cand{\`e}s, Romberg and Tao in \cite{CRT1}, 
Rauhut has shown the following result concerning recovery of sparse 
trigonometric polynomials in \cite[Theorem 2.1 and Section 3.6(b)]{Rau}.

\begin{Theorem}\label{thm:BP}
Let $\Gamma \subset \Z^d$ finite and $T \subset \Gamma$.
Set $M = |T|$ and $D = |\Gamma|$.
Let $X=(x_1,\hdots,x_N)$ be chosen according to one of our two probability models.
If for some $\epsilon > 0$ it holds
\[
N \,\geq\, C M \log (D/\epsilon)
\]
then with probability at least $1-\epsilon$ every trigonometric polynomial
$f \in \Pi_\Gamma(M)$ with Fourier coefficients supported on $T$ 
can be recovered from 
its sample values $f(x_j), j=1,\hdots,N$,
via Basis Pursuit. The constant $C$ is absolute.
\end{Theorem}

In particular, if the sparsity $M$ is small and the dimension 
$D$ large then we may choose
the number $N$ of samples much smaller than $D$ (but larger than $M$), 
and we are still able to recover
a polynomial $f \in \Pi_\Gamma(M)$ exactly -- at least with high probability. 

Note that the theorem is non-uniform in the sense that it does not guarantee
that a single sampling set $\{x_1,\hdots,x_N\}$ is good for all support 
sets $T$ of a certain cardinality $M$. This slight drawback can actually be
removed at the cost of introducing additional $\log$-factors 
by using the concept of the uniform uncertainty principle introduced
in \cite{CRT2,CT}.
For an $N \times D$ matrix $A$ and $M \leq D$ the 
{\it restricted isometry constant}
$\delta_M$ is defined as the smallest number such that for all
subsets $T \subset \{1,\hdots,D\}$ with cardinality $|T| \leq M$
\[
(1-\delta_M) \|x\|_2^2 \leq \|A_T x\|_2^2 \leq (1+\delta_M) \|x\|_2^2
\]
for all coefficient vectors $x$ supported on $T$. Here $A_T$ denotes
the submatrix of $A$ consisting of the columns indexed by $T$.
In \cite{CRT2} the following general recovery theorem for BP was proved.
\begin{Theorem}\label{thm_CRT} Assume that $A$ is a matrix satisfying
\begin{equation}\label{cond_delta}
\delta_M + \delta_{2M} + \delta_{3M} < 1.
\end{equation}
Then every vector $x$ with
at most $M$ non-zero entries can be recovered from $y=Ax$ 
by BP.
\end{Theorem}
 
The above statement is deterministic, but unfortunately, 
it is hard to check
that a deterministic matrix $A$ has small enough $\delta_M$'s. 
So the strategy
is to prove that a random matrix, in particular, our matrix $\F_X$, satisfies
condition (\ref{cond_delta}). Indeed, this was first achieved by Cand{\`e}s
and Tao in \cite{CT} for partial discrete Fourier matrices $\F_X$, and later
improved by Rudelson and Vershynin in \cite{RV06}. Recently, Rauhut \cite{Rau3}
used Rudelson and Vershynin's method to come up with an analog result for
$\F_X$ with samples $X$ drawn from the continuous distribution on
$[0,2\pi]^d$.
More precisely, if
\begin{equation}\label{cond_delta2}
N / \log(N) \geq C_\delta M \log^2(M) \log(D) \log^2(\epsilon^{-1})
\end{equation}
then with probability at least $1-\epsilon$ 
the restricted isometry constant of the $N \times D$ matrix 
$N^{-1/2}\F_X$ satisfies $\delta_M \leq \delta$, 
both for the continuous and discrete probability model, see 
\cite[Theorem 3.2]{Rau3}, \cite{RV06}.
The above condition is particularly good for small $M$, but is always
satisfied if $N \geq C M \log^4(D) \log^2(\epsilon^{-1})$.
(It seems that the $\log^2(\epsilon^{-1})$-factor can be improved
to $\log(\epsilon^{-1})$, see \cite{ruve06-2}).
We note that
Cand{\`e}s and Tao obtained the condition 
$N \geq C M \log^5(D) \log(\epsilon^{-1})$ 
(substitute $\epsilon = cD^{-\rho/\alpha}$ in \cite[Definition 1.12]{CT}
of the uniform uncertainty principle and use \cite[Lemma 4.3]{CT}).  

Combining Theorem \ref{thm_CRT} and condition (\ref{cond_delta2}) gives
a uniform result for recovery of sparse trigonometric polynomials by
BP, see also \cite{CT,Rau3,RV06}. Thus, a single 
sampling set $X$ may be good for all
$f \in \Pi_\Gamma(M)$, in particular for all support sets $T$.

Note that in \cite{CRT2} an extension of Theorem \ref{thm_CRT} 
to noisy and non-sparse situations is provided.

\subsection{Recovery results for greedy algorithms}

Let us first consider recovery by thresholding. Our theorem
reads as follows.
\begin{Theorem}\label{thm:thresh} 
Let $f \in \Pi_\Gamma(M)$ with Fourier coefficients
supported on $T$ with $|T| = M$. 
Define its dynamic range by
\[
R \,:=\, \frac{\max_{k \in T} |c_k|}{\min_{k \in T} |c_k|}.
\]
Choose random sampling points $X=(x_1,\hdots,x_N)$ according to one of our two
probability models.
If for some $\epsilon > 0$
\begin{equation}\label{cond:Thresh}
N \geq C M R^2 \log(4D/\epsilon)
\end{equation}
then with probability at least $1-\epsilon$ 
thresholding recovers the correct support $T$ and, hence, also 
the coefficients $c$.
The constant $C$ is no larger than $17.89$.
\end{Theorem}
Hence, also with thresholding we can recover a sparse trigonometric polynomials
from few random samples. Again the number of samples $N$ is linear in the
sparsity. However, compared to Theorem \ref{thm:BP} for Basis Pursuit 
there is an additional dependence on the dynamic range $R$. Moreover, the above
theorem is non-uniform in contrast to the available results for Basis Pursuit.
This means that recovery is successful with high probability 
only for the given support set $T$ and coefficient vector $c$. It does not 
guarantee that a single sampling set $X$
is good for recovering all sparse trigonometric polynomials.
Indeed, if $N\leq CM^2$ it can be shown that given the random sampling set $X$
there exists with high probability a coefficient vector $c=c(X)$ supported
on $T$ for which thresholding fails \cite{rave07}
(compare also the next section). 

Let us now consider Orthogonal Matching Pursuit. One can expect that
this iterative strategy performs better than the simple Thresholding
algorithm. In particular, the dependence of the number 
of samples $N$
on the dynamic range $R$ as in (\ref{cond:Thresh}) should be removed. 
This is indeed indicated by the next
theorem.

\begin{Theorem}\label{thm:omp} Let $f\in \Pi_\Gamma(M)$ with coefficients
supported on $T$.
Choose random sampling points $X=(x_1,\hdots,x_N)$ according to one of our two
probability models.
If
\[
N \geq C M \log(8D/\epsilon)
\]
then with probability at least $1-\epsilon$ OMP 
selects an element
of the true support $T$ in the first iteration.
The constant $C$ is no larger than $32.62$.
\end{Theorem}
Our numerical experiments indicate that this result should extend
to all iterations so that finally the true sparse polynomial $f$ is recovered.
Unfortunately, we have not yet been able to analyze theoretically 
the further iterations of OMP because of subtle 
stochastic dependency issues.

Again the above theorem is non-uniform and we refer to the next section for
a uniform result, which, however, requires much more samples.
Similarly as for thresholding, one can
show that if $N \leq CM^{3/2}$ then with high probability there exists
an $M$-sparse coefficient vector $c$ depending on the sampling set 
such that OMP fails in the first step \cite{rave07}.
We remark, however, that if OMP selects a wrong element
in some step it might nevertheless 
recover the right polynomial. Indeed, if after
$M' > M$ steps the true support $T$ is contained in the recovered
support $T'$ then the coefficients in $T' \setminus T$ will usually be set to
$0$, since with high probability $\F_{T'X}$ is injective 
(see also Lemma 3.2 in \cite{Rau}). A precise 
analysis when this situation occurs seems to be quite involved, however.

\subsection{Uniform results based on coherence}\label{coherence}

Many results for OMP rely on the so called coherence parameter $\mu$ 
\cite{Tropp_Greed,Tropp_Thesis}. 
It measures the correlation of different columns of the measurement matrix
$\F_X$. It requires that they have unit norm, so we define 
$\widetilde{\phi}_k := N^{-1/2} \phi_k$ resulting in $\|\widetilde{\phi}_k\|_2 = 1$. 
Then $\mu$ is defined as
\begin{equation}\label{def_coh}
\mu \,:=\, \max_{j \neq k} |\langle \widetilde{\phi}_j,\widetilde{\phi}_k\rangle|
\,=\, N^{-1} \max_{j \neq k} |\langle \phi_j, \phi_k\rangle|.
\end{equation}
Reformulating Corollary 3.6 in \cite{Tropp_Greed} for our context yields
the following result. 

\begin{Theorem}\label{thm_coh_rec} 
Assume $(2M-1) \mu < 1$. Then OMP (and also BP) recovers
{\em every} $f \in \Pi_\Gamma(M)$.
\end{Theorem}
\begin{Proof} The matrix $\F_{TX}$ is injective for all $T$ 
with $|T| \leq M$ due to the condition on $\mu$, 
see e.g.~\cite[Proposition 4.3]{Tropp_Thesis}.
This ensures uniqueness of a polynomial $f \in \Pi_\Gamma(M)$ 
having the vector $(f(x_1),\hdots,f(x_N))$ of sample values. 
The rest of the proof 
is the same as the one of Corollary 3.6 in \cite{Tropp_Greed}.
\end{Proof}
A similar result holds also for thresholding, see e.g.~\cite{ICASSP}. 
Again there
is an additional dependence on the ratio of the largest to the smallest
coefficients.

\begin{Theorem}\label{thm:coh:thresh} Assume $(2M-1) \mu < R^{-1}$ for some $R \geq 1$. 
Then thresholding
recovers {\em every} $f \in \Pi_\Gamma(M)$ whose Fourier 
coefficients $c$ satisfy
\[
\frac{\max_{k \in \supp c} |c_k|}{\min_{k \in \supp c} |c_k|} \leq R.
\] 
\end{Theorem}

So clearly, we need to investigate the coherence of the random
matrix $\F_X$. 

\begin{Theorem}\label{thm:coh} Let $X=(x_1,\hdots,x_N)$ be chosen according
to one of our two probability models. 
Suppose that
\begin{equation}\label{cond:coherenceN}
N \geq C (2M-1)^2 \log(4D'/\epsilon),
\end{equation}
where $D' := \#\{j-k: j,k \in \Gamma, j \neq k\} \leq D^2$.
Then with probability at least $1-\epsilon$ 
the coherence of $\F_X$ satisfies 
\[
(2M-1) \mu < 1,
\] 
and consequently OMP recovers {\it every} $M$-sparse trigonometric 
polynomial. The constant satisfies 
$C \leq 4 +\frac{4}{3\sqrt{2}} \approx 4.94$. In case of 
the continuous probability model it can be improved to $C = 4/3$.
\end{Theorem}
Of course, by Theorem \ref{thm:coh:thresh} a similar result
applies also to thresholding. 
Note that always $D \leq D' \leq D^2$. If $\Gamma = \Z_m^d$ and
$X$ is chosen at random from the grid $\frac{2\pi}{m} \Z_m^d$ then
actually $D'=D-1$ because $j-k$ is computed in the ``periodic sense''. 

In contrast to the results of the previous section the above theorem
is uniform in the sense that a single sampling set $X$ is sufficient
to ensure recovery of {\it all} sparse signals. 
However, it clearly has the drawback that now the number of samples $N$
scales quadratically in the sparsity $M$. Apart from perhaps the 
constant $C$ and the logarithmic scaling in the dimension $D$ 
condition (\ref{cond:coherenceN})
actually does not seem to be improvable
if one requires exact recovery of {\it all} sparse
trigonometric polynomials from a single sampling set $X$, see the remark
after Theorem \ref{thm:omp}.
So in this regard there is a crucial difference between 
Basis Pursuit and greedy algorithms. 
For certain applications
a non-uniform recovery result might be enough and then greedy algorithms work
well (and usually much faster than BP), see also the numerical experiments
in Section \ref{sect:num}. 
But if one requires uniformity (and speed is not a concern) 
then Basis Pursuit seems to be the method
of choice.

\subsection{Related Work}
\label{Sec_Related}

{\bf BP with other measurement ensembles.}
Other choices of the measurement matrix $A$ (instead of $\F_X$)
were also investigated by several authors, 
see e.g.~\cite{CT,Donoho1,DTsaig,Donoho2,DT,RV,RV06}. 
For instance, it was shown (see e.g.~\cite{badadewa06,CT,CT2,RV,RV06}) 
that an $N\times D$ random
matrix $A$ with independent Gaussian distributed entries 
(with variance $N^{-1}$) has restricted isometry constant 
$\delta_M \leq \delta$ with probability at least $1-\epsilon$
provided
\begin{equation}\label{Gauss_delta}
N \geq C_\delta M \log\left(\frac{D}{M\epsilon}\right).
\end{equation}
Similar estimates are possible for Bernoulli matrices, i.e., random
matrices with independent $\pm 1$ entries. Condition (\ref{Gauss_delta})
is slightly better than (\ref{cond_delta2}). So for certain
applications of compressed sensing Gaussian / Bernoulli matrices $A$ might 
be useful. However, such ``completely random'' matrices have the disadvantage
of not being structured, hence in contrast to the Fourier case 
no fast algorithms are available for
matrix vector multiplication. Moreover, the ``samples'' $y = Ax$
lack a physical meaning. Of course, this does not matter in 
encoding / decoding problems. However, there are possible 
applications where the samples result as the output of a physical measurement
or an image reconstruction problem \cite{CRT1} and can really
be interpreted as samples of a trigonometric polynomial. Clearly,
in such a case one must use the Fourier matrix $\F_X$.

{\bf Orthogonal Matching Pursuit.} Concerning OMP 
(or greedy algorithms in general) in connection
with compressed sensing much less investigations were
done so far. Gilbert and Tropp proved the following
result \cite[Theorem 2]{Tropp} (set $\epsilon = 2D^{-p}$ to arrive at the formulation below).
\begin{Theorem}\label{Gilbert_Tropp} 
Let $x$ be an $M$-sparse vector in $\R^D$. Let
$A$ be an $N \times D$ random matrix with independent standard Gaussian
entries. Assume that for some $\epsilon > 0$
\[
N \geq 8 M \log(2D/\epsilon).
\]
Then Orthogonal Matching Pursuit reconstructs $x$ from the given data $y = Ax$
with probability exceeding $1-\epsilon$.
\end{Theorem}
The method in \cite{Tropp} actually applies to more general 
types of random matrices
$A$, in particular, to Bernoulli matrices. 
However, it is heavily used
that the columns of $A$ are stochastically independent, 
so unfortunately their approach cannot be applied directly
to the random Fourier matrices $\F_X$. 

The above theorem is non-uniform, i.e.,
recovery is successful with high probability only for a given
sparse vector $x$. It cannot be guaranteed that with high probability 
a single Gaussian 
measurement matrix $A$ suffices for all $M$-sparse vectors $x$.
As in the Fourier case 
Theorem \ref{Gilbert_Tropp} becomes actually false
when requiring reconstruction for all coefficients
supported on a given set $T$, see e.g.~\cite[Section 7]{Donoho2}.

{\bf Ordinary Matching Pursuit}
Ordinary Matching Pursuit is slightly simpler than OMP.
The difference lies in steps 5 and 6
of Algorithm \ref{algo:omp}. Instead of performing the orthogonal
projection over the whole space of previously chosen elements $\phi_{k}$, the
new iterate and residual are computed as
\begin{align*}
d_s &\,=\, d_{s-1} + \langle r_{s-1}, \phi_{k_s} \rangle e_{k_s},\\
r_s &\,=\, r_{s-1} - \langle r_{s-1}, \phi_{k_s} \rangle \phi_{k_s},
\end{align*}
where $e_{k}$ denotes the $k$-th canonical unit vector.
This step is clearly faster than the orthogonal projection step. However,
in contrast to OMP it may now happen that some elements $k_s$ are
selected more than once. So even if MP picks a correct 
element $k \in T$ in every step it usually has not yet reconstructed
the right coefficients after $M = |T|$ steps. 
However, Gribonval
and Vandergheynst showed that exponential
convergence can be guaranteed under the same condition
as for OMP in Theorem \ref{thm_coh_rec}, i.e., $(2M-1) \mu < 1$,
see \cite[Featured Theorems 1 and 2]{grva06}. 
Hence, the uniform recovery result Theorem \ref{thm:coh} has also
an immediate application to ordinary Matching Pursuit.
Moreover, Theorem \ref{thm:omp} is applicable as well, since the first
step of ordinary Matching Pursuit and of Orthogonal Matching Pursuit 
coincide.

{\bf Sublinear Algorithms for Sparse Fourier Analysis.}
Finally, we would like to point out that in the series of papers
\cite{GGIMS,GiMuSt,ZGSD,Zou} an approach from \cite{Mansour} has been
considered.
The algorithms are based on so called isolation and group testing from
theoretical computer science.
It has been proven that one can recover a sparse trigonometric
polynomial $f\in\Pi_{\Gamma}(M)$, $\Gamma=\mathbb{Z}_m$, from $N$ randomly
chosen samples on $\frac{2\pi}{m}\mathbb{Z}_m$ if
\begin{equation*}
  N\ge C M^{\alpha} \log^{\beta}(D) \log^{\gamma}(1/\epsilon)
  \log^{\delta}(R),\qquad \text{with some } \alpha,\, \beta,\, \gamma,\,
  \delta\in\mathbb{N}.
\end{equation*}
These results have been generalized to some extend for multivariate
trigonometric polynomials in \cite{GiMuSt,ZGSD,Zou}.
Moreover, \cite{GiMuSt} indicates that $N\ge C^{\prime} M
\log^{\beta^{\prime}}(D) \log^{\gamma^{\prime}}(1/\epsilon)
\log^{\delta^{\prime}}(R)$ samples suffice if the random sampling set
possesses additional structure.

The striking point in these algorithms is the fact that the computation time
scales ``sublinear'' as $\log^{\beta}(D)$ in the dimension $D$.
However, numerical experiments in \cite{ZGSD,Zou} show even for small
sparsities $M$ a rather large crossover point with the classical FFT which
scales as $D\log(D)$. It seems that in practice these algorithms
require much more samples than OMP and BP, so that the main concern of these
methods is speed rather than a minimal number of samples. 

\section{Proofs}
\label{sect:proofs}

A main tool for our proofs is Bernstein's inequality,
see e.g.~\cite[Lemma 2.2.9]{vawe96} or \cite{be62}.

\begin{Theorem} Let $Y_\ell$, $\ell=1,\hdots,N$, be a sequence of 
independent real-valued random variables with mean zero and variance
$\E [Y_\ell^2] \leq v$ for all $\ell=1,\hdots, N$. Assume that
$|Y_\ell| \leq B$ almost surely. Then
\[
\P\left(|\sum_{\ell=1}^N Y_\ell| \geq x\right) \leq 
2 \exp\left(-\frac{1}{2} \frac{x^2}{Nv + B x/3}\right).  
\]
\end{Theorem}

Bernstein's inequality allows to prove the following
concentration inequality.

\begin{lemma}\label{lem:conc}
Assume that $c$ is a vector supported on $T$. 
Further, assume that the sampling set $X$ is chosen according
to one of our two probability models. 
Then for $j \notin T$ and $x > 0$ it holds
\[
\P\left(|N^{-1} \langle \F_{TX} c, \phi_j \rangle | \geq x\right)
\leq 4\exp\left(- N \frac{x^2}{4\|c\|_2^2 + \frac{4}{3\sqrt{2}}\|c\|_1 x}\right).
\] 
\end{lemma}
\begin{Proof} Note that 
\[
\langle \F_{TX} c, \phi_j \rangle
\,=\, \sum_{k \in T} \sum_{\ell=1}^N c_k e^{i(k-j)\cdot x_\ell} 
\,=\, \sum_{\ell = 1}^N Y_\ell
\]
where
\[
Y_\ell = \sum_{k \in T} c_k e^{i(k-j) \cdot x_\ell},\quad \ell=1,\hdots,N.
\]
Since the $x_\ell$ have the uniform distribution either on the grid
$\frac{2\pi}{m} \Z_m^d$ or on the cube $[0,2\pi]^d$, and since $j \notin T$
it is easy to see that $\E Y_\ell = 0$. Furthermore, 
the random variables $Y_\ell$ are independent and bounded,
\[
|Y_\ell| \leq \sum_{k \in T} |c_k| = \|c\|_1.
\]
Their variance is given by
\[
\E [|Y_\ell|^2] \,=\, \E \left[ \sum_{k,k' \in T} c_k \overline{c_{k'}} e^{i(k-j)\cdot x_\ell} e^{-i(k'-j) \cdot x_\ell}\right]
\,=\, \sum_{k,k' \in T} c_k \overline{c_{k'}} \E [e^{i(k-k')\cdot x_\ell}]
\,=\, \sum_{k \in T} |c_k|^2 \,=\, \|c\|_2^2.
\]
Hereby, we used that $\E[e^{i(k-k')\cdot x_\ell}] = \delta_{k,k'}$. 
Thus, we obtain
\begin{align}
& \P\left(|N^{-1} \langle \F_{TX} c, \phi_\ell \rangle| \geq x\right)
\, =\, \P\left(N^{-1} |\sum_{\ell=1}^N Y_\ell| \geq x\right) \notag\\
&\leq\, \P\left(N^{-1} |\sum_{\ell=1}^N \re (Y_\ell)| \geq x/\sqrt{2}\right) +  
\P\left(N^{-1} |\sum_{\ell=1}^N \im (Y_\ell)| \geq x/\sqrt{2}\right).\notag
\end{align}
Applying Bernstein's inequality to $\re(Y_\ell)$ and $\im(Y_\ell)$, 
$\ell = 1,\hdots,N$, (using that $|\re(Y_\ell)|\leq |Y_\ell|$ and 
$|\im(Y_\ell)| \leq |Y_\ell|$) we obtain
\begin{align}
\P(|N^{-1} \langle \F_{TX} c, \phi_\ell \rangle| \geq x)
\,&\leq\, 4 \exp\left(-\frac{1}{2} \frac{N^2x^2/2}{N\|c\|_2^2 +\frac{1}{3\sqrt{2}} \|c\|_1 N x}\right)\notag\\
&=\, 4 \exp\left(-N \frac{x^2}{4 \|c\|_2^2 + \frac{4}{3\sqrt{2}} \|c\|_1 x}\right).\notag
\end{align}
This concludes the proof.
\end{Proof}

\subsection{Proof of Theorem \ref{thm:thresh}}

Thresholding recovers the correct support if
\[
\min_{j \in T} |N^{-1} \langle \phi_j, \F_{TX} c\rangle| > 
\max_{k \notin T} |N^{-1} \langle \phi_k, \F_{TX} c\rangle|.
\]
Observe that if $l\in T$ then the triangle inequality yields
\begin{align*}
  |N^{-1} \langle \phi_l, \F_{TX} c\rangle| &= |c_l + N^{-1}\langle \phi_l,
  \F_{(T\setminus\{l\})X} c_{T\setminus \{l\}}\rangle|\\
  &\ge
  \min_{j \in T}|c_j| - \max_{j \in T} |N^{-1}\langle \phi_j,
  \F_{(T\setminus\{j\})X} c_{T\setminus \{j\}}\rangle|,
\end{align*}
where $c_{T \setminus \{l\}}$ denotes the vector $c$ restricted to the indices
in $T \setminus \{l\}$.
Thus, thresholding is certainly successful if
\[
\max_{j \in T} |N^{-1} \langle \phi_j, \F_{(T \setminus \{j\})X} c_{T\setminus \{j\}}\rangle| <
\min_{j \in T} |c_j|/2 
\quad\text{and}\quad
\max_{k \notin T} |N^{-1} \langle \phi_k, \F_{TX} c \rangle| <
\min_{j \in T} |c_j|/2.
\]
We conclude that the probability that thresholding is {\em not} successful
can be upper bounded by
\begin{align}
\P&\left(\max_{j \in T} |N^{-1}\langle \phi_j, \F_{(T \setminus \{j\})X} c_{T \setminus \{j\}}\rangle| \ge 
\min_{j \in T} |c_j|/2\right)
+ \P\left(\max_{k \notin T} |N^{-1} \langle \phi_k, \F_{TX} c \rangle| \ge
\min_{j \in T} |c_j|/2 \right)\notag\\
&\leq \sum_{j \in T} \P\left(|N^{-1}\langle \phi_j, \F_{(T \setminus \{j\})X} c_{T \setminus \{j\}}\rangle| \ge
\min_{j \in T} |c_j|/2 \right)
+ \sum_{k \notin T} \P\left( |N^{-1} \langle \phi_k, \F_{TX} c \rangle| \ge
\min_{j \in T} |c_j|/2 \right). \notag
\end{align}
Applying Lemma \ref{lem:conc} we can bound the probability that
thresholding fails by
\begin{align}
&4|T| \exp\left(-N \frac{\min |c_j|^2/4}{4 \max\|c_{T
      \setminus \{j\}}\|_2^2 + \frac{4}{3\sqrt{2}} \max \|c_{T\setminus \{j\}}\|_1 \min |c_j|/2}\right) \notag\\
& + 4(D-|T|) \exp\left(-N \frac{\min |c_j|^2/4}{4 \|c\|_2^2 + \frac{4}{3\sqrt{2}} \|c\|_1\min |c_j|/2}\right)\notag\\
\leq &  4 D \exp\left(-N \frac{\min |c_j|^2}{16 \|c\|_2^2 + \frac{8}{3\sqrt{2}} \|c\|_1 \min |c_j|}\right) \notag\\
\leq & 4 D \exp\left(-N \frac{ \min |c_j|^2}{16 M \max |c_j|^2 + \frac{8}{3\sqrt{2}}
M \max |c_j| \min |c_j|}\right) \label{prob:estim}\\
\leq & 4D \exp\left(-\frac{N}{M} \frac{\min |c_j|^2}{\max |c_j|^2} \frac{1}{16 + \frac{8}{3\sqrt{2}}}\right), \label{prob:estim2}
\end{align}
where $\min|c_j|$ and $\max|c_j|$ are always taken over $j\in T$.
We note that in (\ref{prob:estim}) the following obvious estimates
were used,
\[
\|c\|_2^2 \leq M \max_{j \in T} |c_j|^2,\quad \mbox{and}\quad
\|c\|_1 \leq M \max_{j \in T} |c_j|.
\]
Requiring that the term in (\ref{prob:estim2}) is less than $\epsilon$
and solving for $N$ shows the claim.

\subsection{Proof of Theorem \ref{thm:omp}}

As basic ingredient of the proof we will use that  
a submatrix $\F_{TX}$ is well-conditioned under
mild conditions on the number of sampling points $N$ and
sparsity $M$. The corresponding result is taken from 
\cite{grrapo06}. It is based on an analysis in \cite{Rau}
of the expectation of the Frobenius norm of 
high powers of $N I - \F_{TX}^* \F_{TX}$, which uses a combinatorial
argument based on set partitions and is inspired by \cite{CRT1}.

\begin{Theorem}\label{thm:eigvals} 
Let $T\subset\Z^d$ be of size $|T| = M$ and let $x_1,\hdots,x_N$ be chosen
according to one of our two probability models.
Choose $\epsilon, \delta \in (0,1)$ and assume
\[
\left \lfloor \frac{\delta^2 N}{3e M} \right\rfloor \,\geq\,
\log (c(\delta) M/\epsilon),
\]
where $c(\delta) = (1-e^{-1}\delta^2)^{-1} \leq (1-e^{-1})^{-1} \approx 1.582$.
Then with probability at least $1-\epsilon$ the minimal and
maximal eigenvalue
of $\F_{TX}^* \F_{TX}$ satisfy 
\begin{equation}\label{eigs}
1-\delta \leq \lambda_{\min}(N^{-1} \F_{TX}^* \F_{TX}),
\quad \mbox{and} \quad \lambda_{\max}(N^{-1} \F_{TX}^* \F_{TX})
\leq 1+\delta.
\end{equation}
\end{Theorem}
Written in a more compact form, if
\[
N \geq C \delta^{-2} M \log(M/\epsilon)
\]
then (\ref{eigs}) holds with probability at least $1-\epsilon$. 
It seems that also techniques developed by Vershynin \cite{Versh} can be used 
to provide a version of Theorem \ref{thm:eigvals}.

Now let us turn to the proof of Theorem \ref{thm:omp}.
(Orthogonal) Matching Pursuit 
recovers an element of the support $T$ in the
first step if
\[
\frac{ \max_{j \notin T} |\langle \phi_j, \F_{TX} c\rangle|}{
\max_{k \in T} |\langle \phi_k,\F_{TX} c \rangle|} \,=\, 
\frac{\max_{j \notin T}|N^{-1}\langle \phi_j, \F_{TX} c\rangle|}
{\|N^{-1}\F_{TX}^* \F_{TX} c \|_\infty} < 1.
\]
Assume for the moment that we are on the event
that the minimal eigenvalue of $N^{-1}\F_{TX}^* \F_{TX}$ is bounded
from below by $1-\delta$, cf. Theorem \ref{thm:eigvals}.
Then
\[
\|N^{-1} \F_{TX}^* \F_{TX} c\|_\infty \geq M^{-1/2}
\|N^{-1} \F_{TX}^* \F_{TX} c\|_2 \geq M^{-1/2}(1-\delta) \|c\|_2.
\]
Hence, we need to bound the probability that
\[
\max_{j \notin T} |N^{-1} \langle \phi_j, \F_{TX} c\rangle|
< \frac{1-\delta}{\sqrt{M}} \|c\|_2.
\]
Using Lemma \ref{lem:conc} we can estimate the probability
of the complementary event by
\begin{align}
&\P\left(\max_{j \notin T} |N^{-1} \langle \phi_j, \F_{TX} c\rangle|
\geq \frac{1-\delta}{\sqrt{M}} \|c\|_2\right)\notag\\
&\leq \sum_{j \notin T} \P\left(|N^{-1} \langle \phi_j, \F_{TX} c\rangle|
\ge \frac{1-\delta}{\sqrt{M}} \|c\|_2\right)\notag\\
&\leq 4 (D-M) \exp\left(-\frac{N}{M}
\frac{(1- \delta)^2 \|c\|_2^2}{4\|c\|_2^2 + 
\frac{4}{3\sqrt{2}}\|c\|_1\|c\|_2(1-\delta)/\sqrt{M}}\right) \notag\\
& \leq 4 (D-M) \exp\left(-\frac{N}{M}
\frac{(1-\delta)^2}{4 + \frac{4}{3\sqrt{2}}(1-\delta)}\right).\notag
\end{align}
In the last step, we applied the Cauchy-Schwarz inequality $\|c\|_1 \leq
\sqrt{M} \|c\|_2$.

Altogether, the probability that OMP fails in the first step
can be bounded by
\begin{equation}\label{prob:estim3}
\P(\lambda_{\min}(N^{-1} \F_{TX}^* \F_{TX}) < 1-\delta)
+ 4 (D-M) \exp\left(-\frac{N}{M}
\frac{(1-\delta)^2}{4 + \frac{4}{3\sqrt{2}}(1-\delta)}\right).
\end{equation}
Now choose $\delta = 1/2$. Then by Theorem \ref{thm:eigvals} 
the first term above can be bounded by
$\epsilon/2$ provided
\begin{equation}\label{cond:1}
\left\lfloor \frac{N}{12e M} \right\rfloor
\geq \log\left(\frac{2}{1-1/(4e)} M/\epsilon\right). 
\end{equation}
Further, requiring that also the second term in (\ref{prob:estim3})
be less than $\epsilon/2$ yields
\begin{equation}\label{cond:2}
N \geq \left(16 + \frac{8}{3\sqrt{2}}\right)M \log(8(D-M)/\epsilon). 
\end{equation}
Since always $M\leq D$ there exists a constant $C\le 32.62$ such that
\[
N \geq C M \log(8D/\epsilon)
\]
implies that both (\ref{cond:1}) and (\ref{cond:2}) are satisfied.
This concludes the proof.

\begin{remark}\label{rem:iterations}
Clearly, we would like to analyze also
the further steps of OMP. However, starting with the second
step the coefficients $c^{(s)}$ (as well as the current support set $T_s$) 
of the current residual $r_s$
depend on the chosen random sampling points. In this situation
Lemma \ref{lem:conc} cannot be applied directly anymore and
the analysis becomes more difficult.
\end{remark}

\subsection{Proof of Theorem \ref{thm:coh}}

Let $k,j \in \Gamma$ with $k \neq j$. Observe that
\[
N^{-1} \langle \phi_j,\phi_k\rangle \,=\, N^{-1} \sum_{\ell=1}^N
e^{i(j-k)\cdot x_\ell} \,=\, N^{-1} \sum_{\ell=1}^N Y_\ell^{(j-k)}, \qquad
Y_\ell^{(j-k)} := e^{i(j-k)\cdot x_\ell}.
\]
The random variables $Y_\ell^{(j-k)}$ obey $\E[Y_\ell^{(j-k)}] = 0$,
$|Y_\ell^{(j-k)}| = 1$, and $\E[|Y_\ell^{(j-k)}|^2] = 1$.
Hence, the union bound and Bernstein's inequality applied to
$\re(Y_\ell^{(j-k)})$ and $\im(Y_\ell^{(j-k)})$, see also the proof of Lemma
\ref{lem:conc}, yields
\begin{equation}\label{derive:Bernstein}
  \P(\max_{\genfrac{}{}{0pt}{}{j,k\in\Gamma}{j\ne k}} N^{-1} |\sum_{\ell=1}^N Y^{(j-k)}_\ell| \ge x)
  \le \sum_{\genfrac{}{}{0pt}{}{j,k\in\Gamma}{j\ne k}} \P(N^{-1} |\sum_{\ell=1}^N Y^{(j-k)}_\ell| \ge x)
  \le 4D' \exp\left(\frac{-N x^2}{4+\frac{4}{3\sqrt{2}}x}\right).
\end{equation}
Setting $x = (2M-1)^{-1}$ yields
\[
\P((2M-1)\mu \ge 1) \leq 4D' \exp\left(-\frac{N}{(2M-1)^{2}} \frac{1}{4+\frac{4}{3\sqrt{2}}}\right).
\]
Requiring that the latter term be less than $\epsilon$ and solving
for $N$ shows the claim.

In case of the continuous probability model (1) it is possible to improve the
constant.
In \cite{pe93} a sharp moment bound (Khintchine's inequality) for the sum of
the variables $Y_\ell^{(j-k)}$ is provided, i.e.,
\begin{equation}\label{Khintchine}
\E[|\sum_{\ell=1}^N Y_{\ell}^{(j-k)}|^{2p}] \leq p! N^p.
\end{equation}
Markov's inequality yields for $\kappa \in (0,1)$, see also \cite[Proposition 15]{tro-random}, 
\begin{align}
&\P(N^{-1} |\sum_{\ell=1}^N Y_\ell^{(j-k)}| \ge x) \leq 
e^{-\kappa N x^2} \E \left[\exp\left( \kappa N^{-1} |\sum_{\ell=1}^N Y_\ell^{(j-k)}|^2\right)\right]\notag\\
&=\, e^{-\kappa N x^2} \sum_{p=0}^\infty \kappa^p N^{-p} 
\frac{\E[ |\sum_{\ell=1}^N Y_\ell^{(j-k)}|^{2p}]}{p!}
\,\leq\,e^{-\kappa N x^2} \sum_{p=0}^\infty \kappa^p 
\,=\, \frac{e^{-\kappa N x^2}}{1-\kappa}.
\notag
\end{align}
Using this estimate instead of Bernstein's inequality, the above derivation
\eqref{derive:Bernstein} yields
\[
\P((2M-1)\mu \ge 1) \leq \frac{1}{1-\kappa} D' \exp\left(-\frac{N}{(2M-1)^{2}} \kappa\right).
\]
Choosing for instance $\kappa = 3/4$ results in the slightly improved
condition, cf. (\ref{cond:coherenceN}),
\[
N \geq \frac{4}{3}(2M-1)^2 \log(4D'/\epsilon).
\]

If the sampling points are chosen at random from the grid 
$\frac{2\pi}{m}\Z_m^d$ then the moment bound (\ref{Khintchine})
remains valid only for $p=1,\hdots,m-1$. One still expects
that also in this case the constant $C$ can be improved 
(possibly depending on $m$), 
but the method above is not applicable directly any more.

%-----------------------------------------------------------------------------
\section{Implementation and Numerical Experiments}\label{sect:num}

\subsection{Analysis of the time complexity}
Orthogonal Matching Pursuit, cf.~Algorithm \ref{algo:omp}, contains two costly
computations.
Step 4 multiplies the adjoint measurement matrix $\F_X^*$ with the current
residual vector $r_s$.
When drawing the sampling set from the lattice $\frac{2\pi}{m}\Z^d_m$ or from
the cube $[0,2\pi]^d$, we use a zero padded fast Fourier transform (FFT)
or the nonequispaced FFT \cite{DPT,kupo02C}, respectively.
In both cases, the total costs of this step in one iteration is ${\cal
  O}(D\log D)$.
Note furthermore, that if the maximum in step 4 occurs  at several indices the
algorithm chooses one of them.

Step 5 solves in each iteration a least squares problem
\begin{equation*}
  \left\|\F_{T_s X} d_s - \ff\right\|_2\stackrel{d_s}{\rightarrow}\min.
\end{equation*}
A straightforward implementation yields costs ${\cal O}(M N^2)$ per
iteration.
We accelerated this step by a QR factorization of $\F_{T_s X}$ computed from
the factorization of $\F_{T_{s-1} X}$, cf.~\cite[pp. 132]{Bj96}, which reduces
the costs to ${\cal O}(N^2)$ per iteration
Alternatively, we use the iterative algorithm LSQR \cite{PaSa82} which solves
the least squares problem with only ${\cal O}(MN)$ floating point operations
since the matrices $\F_{T_s X}$ have uniformly bounded condition numbers,
cf. Theorem \ref{thm:eigvals}.

Clearly, if OMP succeeds, the algorithm takes $M$ outer iterations.
Two reasonable choices for stopping criteria are a maximum number of
iterations (assuming an upper bound on the sparsity $M$ is known) or a
residual tolerance $\varepsilon$ (or a combination of both).
In any case the algorithm will do no more than $N$ iterations.
Assuming $M$ outer iterations, a reasonable sparsity $M={\cal O}(\sqrt{D})$, and
$N={\cal O}(M\log D)$ sampling points as suggested by Theorem \ref{thm:omp}, 
Algorithm \ref{algo:omp} (using LSQR) has a 
total cost of ${\cal O}(D^{1.5} \log D)$ arithmetic operations.

Thresholding, cf.~Algorithm \ref{algo:thresh}, avoids the iterative
procedure and thus takes ${\cal O}(D \log D)$ operations if we assume again
$M={\cal O}(\sqrt{D})$ and $N={\cal O}(M\log D)$, see also Theorem
\ref{thm:thresh}.
Note that the computation of the inner products $\F_{X}^* r_s$ becomes the
dominant task in both schemes if $M=o(\sqrt{D})$.

\subsection{Matlab toolbox for Thresholding, OMP, and BP}
For the reader's convenience, we provide an efficient and reliable
implementation of the presented algorithms for the univariate case in
Matlab.
Following the common accepted concept of {\em reproducible research}, all
numerical experiments are included in our publicly available toolbox
\cite{KuRa_OMP}.
The toolbox comes with a simple version of the nonequispaced FFT.

All examples testing the Basis Pursuit principle use the optimization tools of
CVX \cite{GrBo_CVX}, L1MAGIC \cite{roca_l1MAGIC}, or MOSEK
\cite{MOSEK}, respectively.
If the vector of Fourier coefficients is assumed to be real valued, the
$\ell^1$-minimization problem is reformulated as a linear program, whereas for
complex valued coefficients the corresponding second order cone problem is set
up.
CVX and MOSEK handle both problems but the matrix $\F_X$ has to be stored
explicitly and no fast matrix vector multiplications by FFTs can be employed.
In contrast, L1MAGIC includes the use of function handles to avoid this memory
bottleneck and reduces the computation time from ${\cal O}(DN)$ to
${\cal O}(D\log D)$ when multiplying with the matrix $\F_X$.
Unfortunately, the solver for equality constraint $\ell^1$-minimization of
this package supports only real valued coefficients.

\subsection{Testbed and examples}
All numerical results were obtained on a Intel PentiumM with 1.6GHz,
512MByte RAM running OpenSUSE Linux kernel 2.6.13-15-default and
MatLab7.1.0.183 (R14) Service Pack 3.
Subsequently, we compare our implementation of Algorithm \ref{algo:omp} and
Algorithm \ref{algo:thresh} with different Basis Pursuit implementations for
the univariate case.

For given dimension $D\in2\N$ we set $\Gamma=\{-D/2,-D/2+1,\hdots,D/2-1\}$.
The support  $T\subset\Gamma$ of the Fourier coefficients is then chosen
uniformly at random among all subsets of $\Gamma$ of size $M$.
We use two choices of (pseudo-)random Fourier coefficients.
If not stated otherwise, the real as well as the imaginary part of the
supported coefficients is drawn from a normal distribution with mean zero and
standard deviation one.
For some experiments with Thresholding, we choose coefficients with modulus
one and draw the phase which from the uniform distribution on $[0,2\pi]$.

The sampling points $x_j$ are drawn uniformly from the interval $[0,2\pi]$ for
the continuous probability model, denoted by NFFT subsequently.
Within the discrete probability model a subset of size $N$ is chosen uniformly
among all subsets of $\{0,\frac{2\pi}{D},\hdots,\frac{2\pi(D-1)}{D}\}$ with
size $N$, denoted by FFT.

\begin{figure}%[ht!]
  \centering
  \subfigure[Success rate vs.~sparsity $M=1,\hdots,40$ for $D=100$, $N=40$,
  $|c_k|=1$ for $k\in T$, FFT.]   
  {\includegraphics[width=0.45\textwidth]{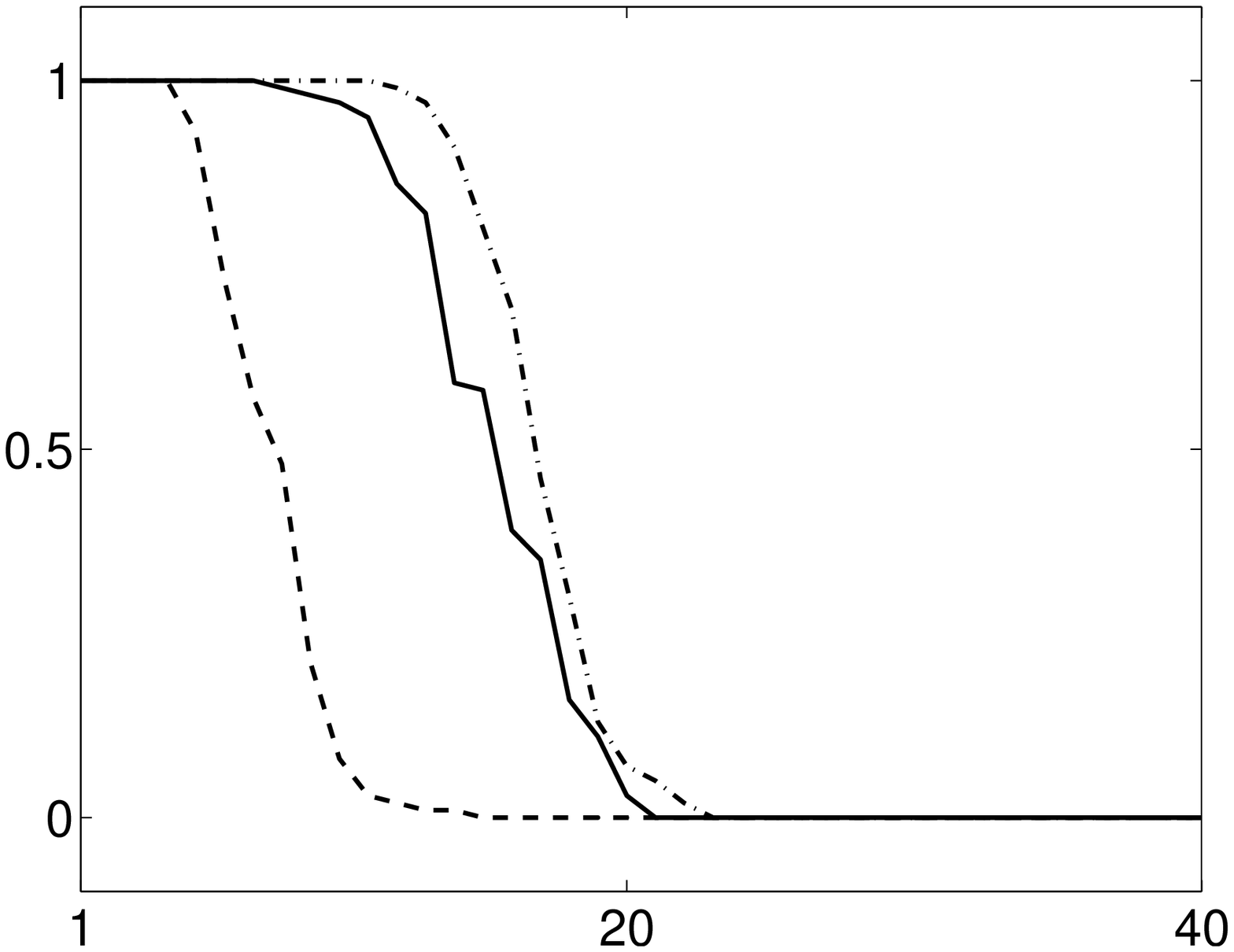}}\hfill
  \subfigure[Success rate vs.~sparsity $M=1,\hdots,40$ for $D=100$, $N=40$, FFT.]   
  {\includegraphics[width=0.45\textwidth]{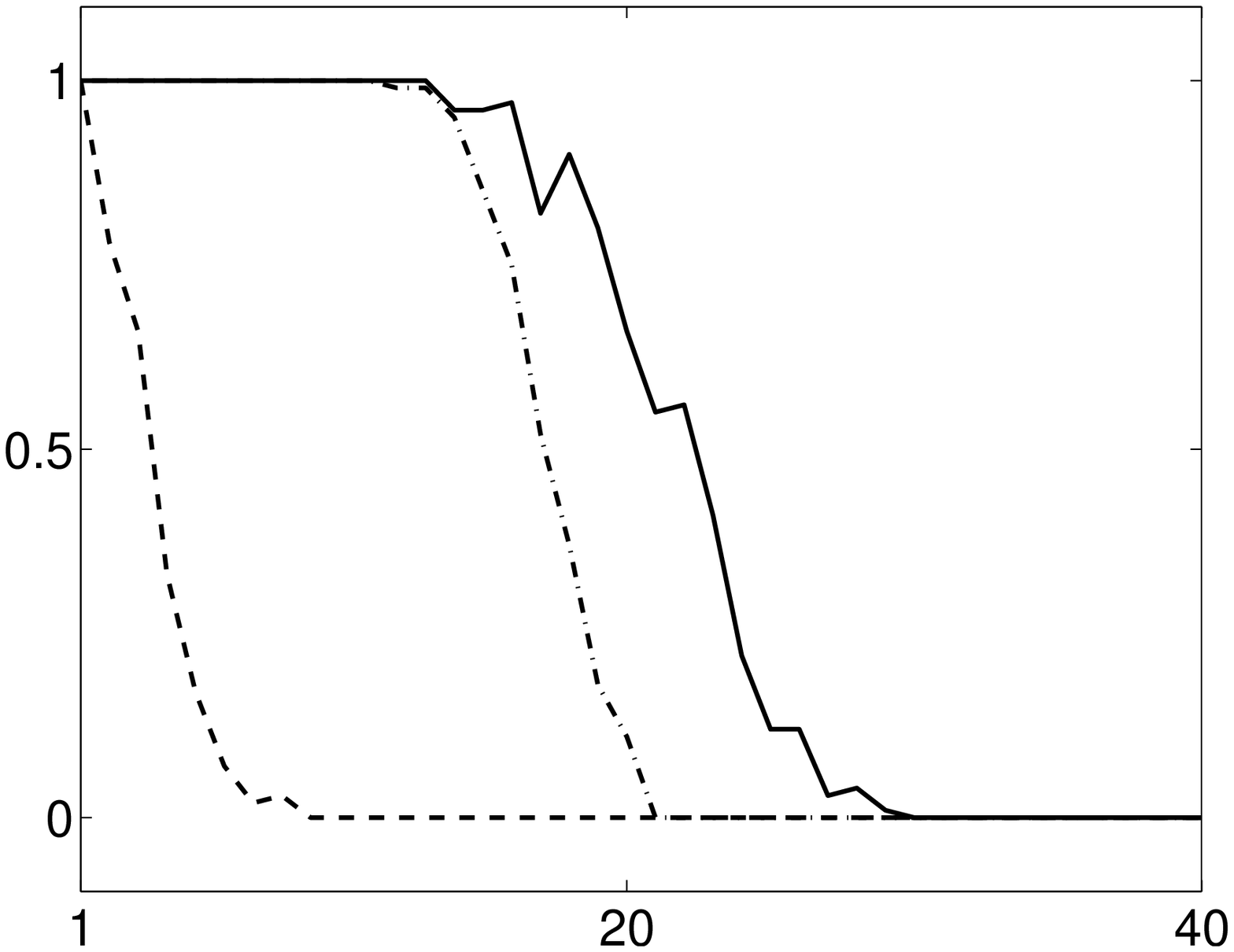}}\\
  \subfigure[Success rate vs.~sparsity $M=1,\hdots,40$ for $D=100$, $N=40$, NFFT.]  
  {\includegraphics[width=0.45\textwidth]{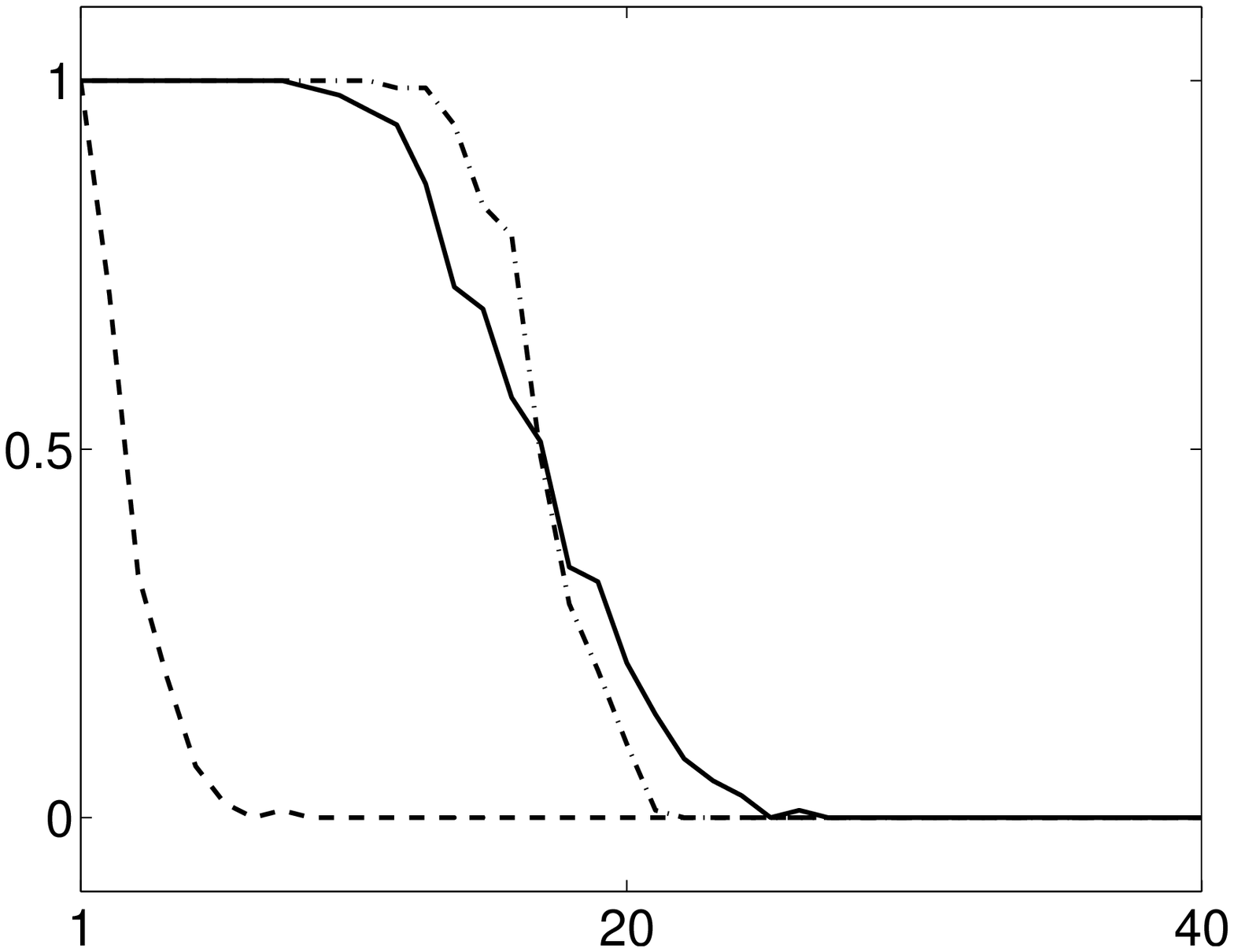}}\hfill
  \subfigure[Computation time in seconds vs.~sparsity $M=1,\hdots,40$ for $D=100$, $N=40$, NFFT.]
  {\includegraphics[width=0.45\textwidth]{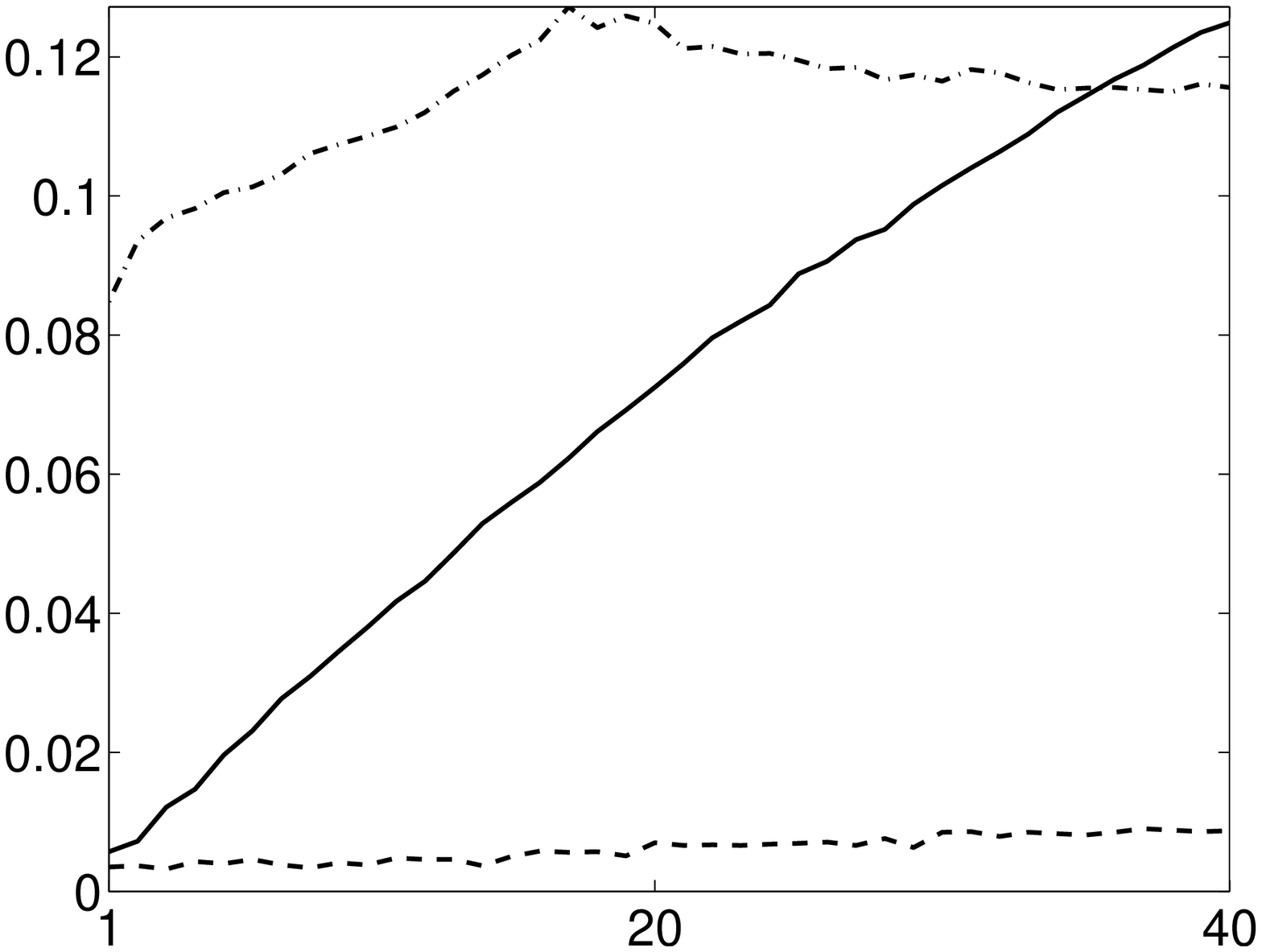}}\\
  \subfigure[Success rate vs.~sparsity $M=1,\hdots,40$ for $D=1000$, $N=80$, FFT.]   
  {\includegraphics[width=0.45\textwidth]{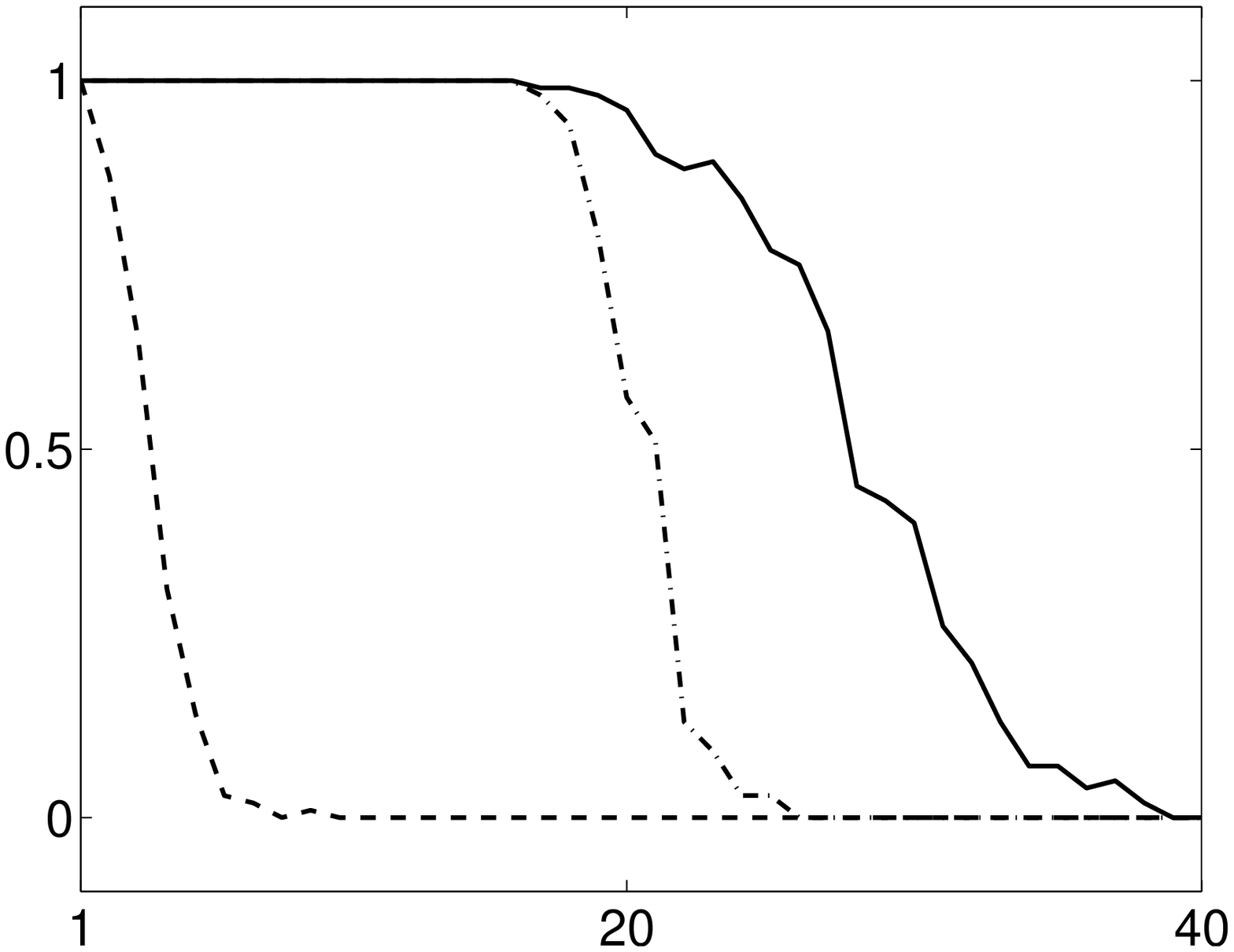}}\hfill
  \subfigure[Computation time in seconds vs.~sparsity $M=1,\hdots,40$ for $D=1000$, $N=80$, FFT.]
  {\includegraphics[width=0.45\textwidth]{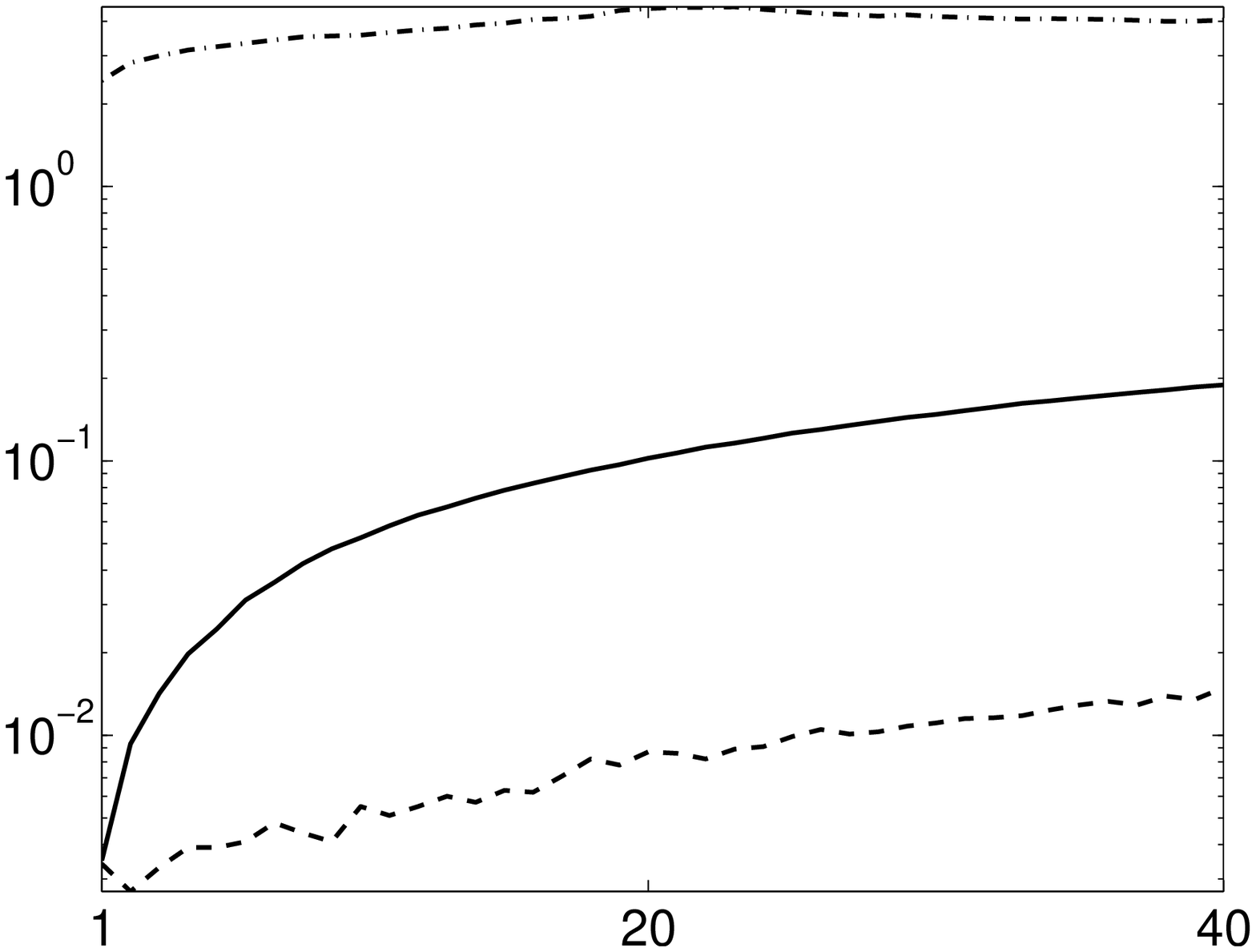}}\\
  \caption{Success rate and computation time of Thresholding (dashed),
    OMP (solid), and BP (dash-dot) with respect to an increasing number $M$ of
    non-zero Fourier coefficients.
    The number of samples $N$ and the dimension $D$ remain fixed; $100$
    runs have been conducted for each setting.\label{fig:ex2}}
\end{figure}

\begin{example}
% File    example2.m
In our first example, we compare the ability of OMP, Thresholding, and BP to
reconstruct sparse trigonometric polynomials with complex valued coefficients
for the dimension $D=50$.
We draw a set $T$ of size $M\in\{1,2,\hdots,40\}$ and $M$ complex valued
Fourier coefficients (with modulus one and a uniformly distributed phase in
Figure \ref{fig:ex2}(a)).
Furthermore, we choose $N=40$ sampling points within the discrete or the
continuous probability model.
The samples $(x_j,f(x_j)),\,j=1,\hdots,N,$ of the corresponding trigonometric
polynomial and the dimension $D$ are the input for OMP, Thresholding, and BP.
OMP and Thresholding use (updated) QR factorizations to solve the least
squares problems, BP uses the MOSEK-package \cite{MOSEK} to solve the second
order cone problem.
The output $d_k\in\C,\,k\in \Gamma,$ of all algorithms is compared to the
original vector of Fourier coefficients.
Repeating the experiment $100$ times for each number $M$ of non-zeros, we
count how often each algorithm is able to reconstruct the given coefficients.
Furthermore, the average CPU-time used by each algorithm with respect
to the number of non-zero coefficients is shown.
The same experiment is done for the dimension $D=1000$.

As readily can be seen from Figure \ref{fig:ex2}(d) and \ref{fig:ex2}(f),
Thresholding, OMP, and BP differ by orders of magnitude in their computation
times.
Note furthermore, that the number of outer iterations and hence the
computation time growths linearly with the number of non-zeros for OMP.
Regarding the rate of successful reconstructions, we observe the following:
Thresholding fails already for a moderate number of non-zero coefficients,
even if we use Fourier coefficients with modulus one, i.e., a dynamic range
$R=1$, cf.~Figure \ref{fig:ex2}(a).
Basis Pursuit and Orthogonal Matching Pursuit achieve much better
reconstruction rates, where surprisingly OMP outperforms BP for larger
dimensions $D$, cf.~Figure \ref{fig:ex2}(e).
Moreover, BP shows the same behavior for different dynamic ranges, the
discrete and the continuous probability model,  cf. Figure
\ref{fig:ex2}(a,b,c).
In contrast, OMP seems to take advantage of a larger dynamic range of the
coefficients in the FFT-situation, cf. Figure \ref{fig:ex2}(b).
\end{example}

\begin{figure}[ht!]
  \centering
  \subfigure[Success rate of BP vs.~sparsity $M=1,\hdots,40$ for $D=100$, $N=30$.]   
  {\includegraphics[width=0.45\textwidth]{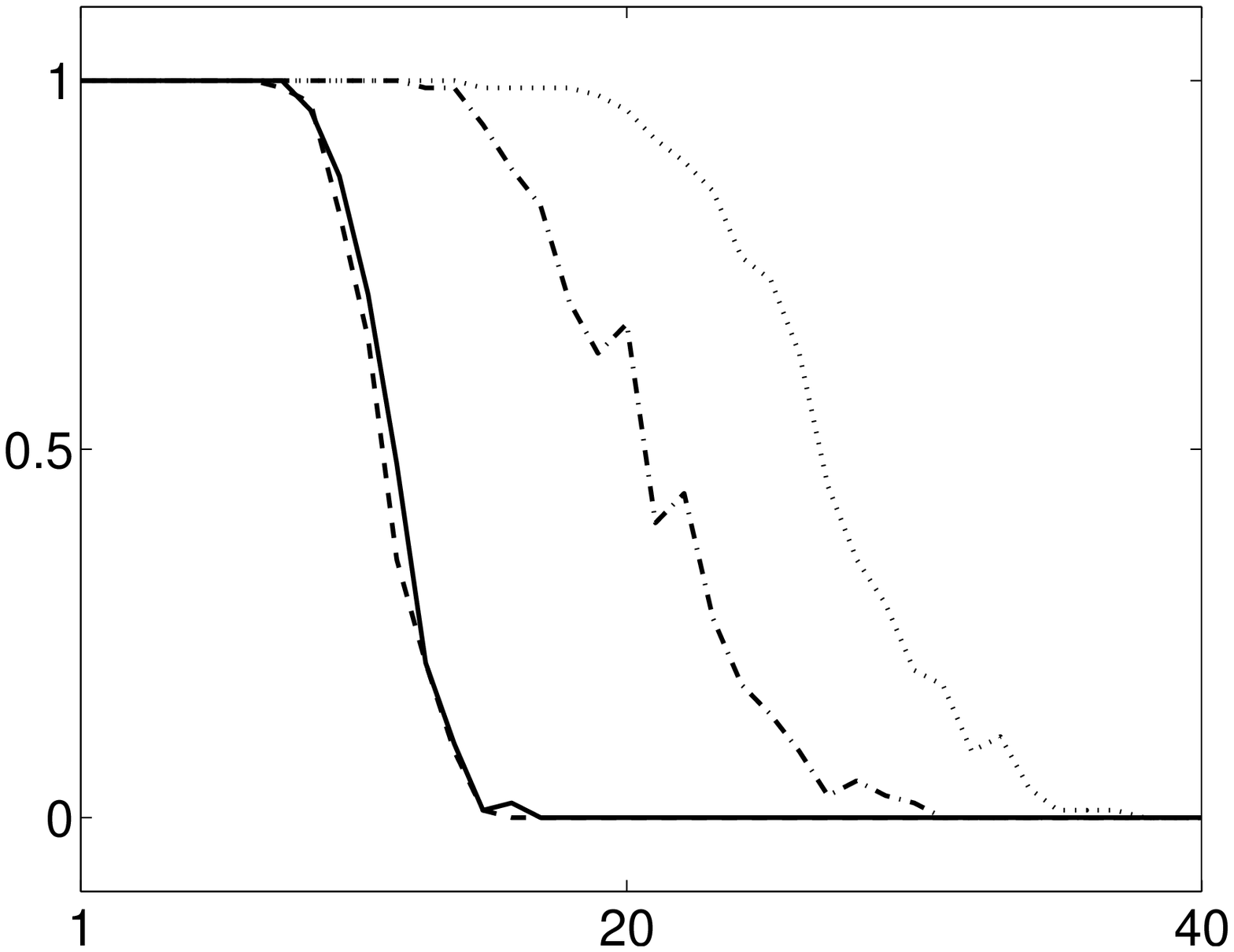}}\hfill
  \subfigure[Computation time of BP vs.~sparsity $M=1,\hdots,40$ for $D=100$, $N=30$.]
  {\includegraphics[width=0.45\textwidth]{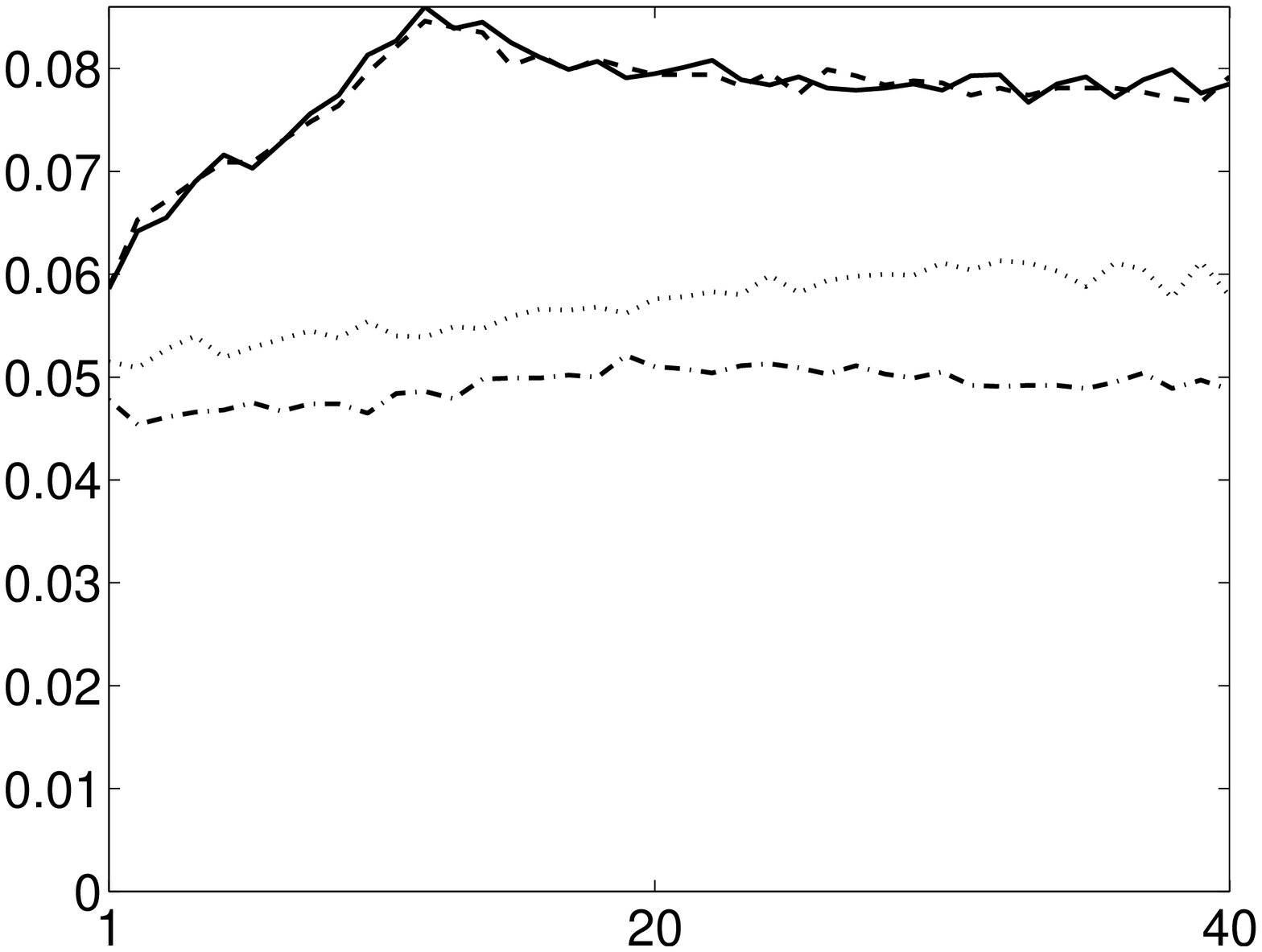}}\\
  \caption{Success rate and computation time of Basis Pursuit for
    complex-valued coefficients with FFT (solid) and NFFT (dash-dot) and
    real-valued coefficients with FFT (dashed) and NFFT (dotted), with
    respect to an increasing number $M$ of non-zero coefficients.
    The number of samples $N$ and the dimension $D$ remain fixed; $100$
    runs have been conducted for each setting.\label{fig:ex3}}
\end{figure}

\goodbreak
\begin{example}
% File    example3.m
The purpose of the second experiment is the comparison for real-valued
and complex-valued coefficients.
We consider the Basis Pursuit principle exclusively.
As in the previous example, we let $D=100$ but use only $N=30$ sampling
points.
Again, we repeat the experiment $100$ times for each level of sparsity.
The results in Figure \ref{fig:ex3} reveal the following:
Besides the easier implementation and a speed up by $25$ percent, the
additional assumption to have real valued coefficients indeed saves roughly
half of the needed samples to recover a sparse trigonometric polynomial.
\end{example}

\begin{figure}[ht!]
  \centering
  \subfigure[Success rate of OMP vs.~sparsity $M=1,\hdots,40$ for $D=1024$ and
  $N=\theta M$ with $\theta=2.5$ (solid), $\theta=3$ (dashed), and $\theta=3.5$
  (dash-dot).]
  {\includegraphics[width=0.45\textwidth]{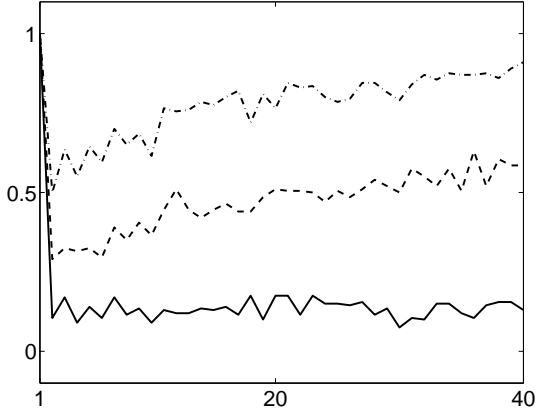}}\hfill
  \subfigure[Oversampling factor $\theta$ vs.~dimension $D$ for Thresholding
  ($|c_k|=1$ for $k\in T$), $90\%$ success.]
  {\includegraphics[width=0.45\textwidth]{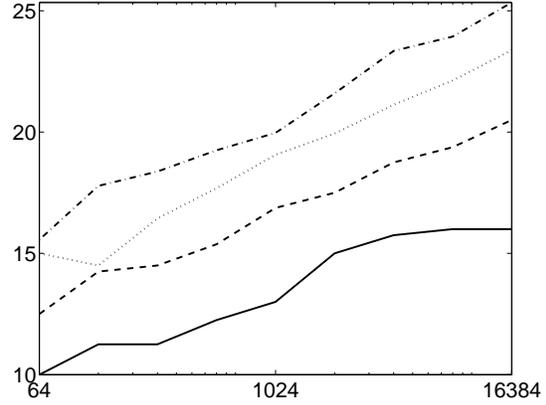}}\\
  \subfigure[Oversampling factor $\theta$ vs.~dimension $D$ for Orthogonal
  Matching Pursuit, $90\%$ success.]
  {\includegraphics[width=0.45\textwidth]{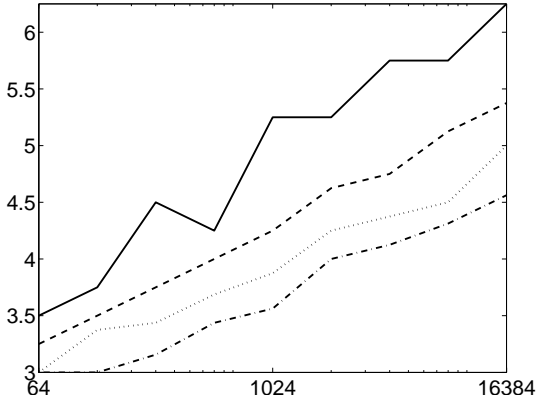}}\hfill
  \subfigure[Oversampling factor $\theta$ vs.~dimension $D$ for Basis Pursuit,
  $90\%$ success.]
  {\includegraphics[width=0.45\textwidth]{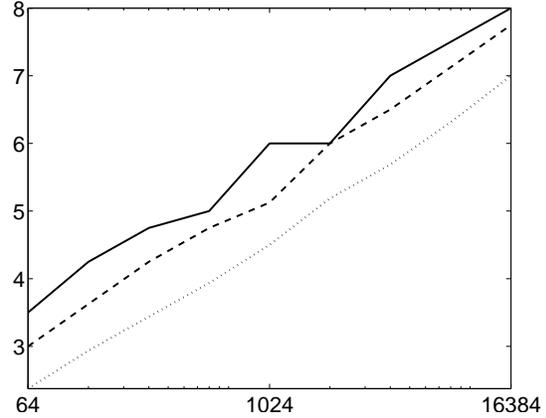}}\\
  \caption{Relation between the generalized oversampling factor $\theta=N/M$
    and the success rate for the continuous probability model.
    Figures (b)-(d): Oversampling factors $\theta$ necessary to reach a
    success rate  of $90$ percent with respect to the dimension $D$ and fixed
    numbers $M=4$ (solid), $M=8$ (dashed), $M=16$ (dotted), and $M=32$
    (dash-dot) of non-zero coefficients. \label{fig:ex4}}
\end{figure}

\begin{example}
% File example6.m  example5_th.m example5_omp.m example5_bp.m 
This example comments on the generalized oversampling factor $\theta=N/M$ and
its relation to the rate of successful reconstructions.
For a fixed dimension $D=1024$ and a fixed generalized oversampling factor
$\theta>0$, we randomly draw a support set $T\subset\Gamma$ of increasing
sizes $M=1,\hdots,40$, random Fourier coefficients (with modulus one and a
uniformly distributed phase for Thresholding), and sampling sets in the
continuous probability model of size $N=\theta M$.
For $200$ runs of each experiment, we count the number of perfect
reconstructions by OMP after exactly $M$ steps.
As Figure \ref{fig:ex4}(a) reveals, the success rate stays (almost) constant
or might even increase slightly for an increasing number of non-zero
coefficients if the ratio $\theta=N/M$ remains constant, cf. Theorem
\ref{thm:omp} and \ref{thm_coh_rec}.

In the second part of this example, we test Thresholding, OMP, and BP.
Our concern is the dependence of the ratio $\theta=N/M$ to reach a certain
success rate when the dimension $D$ varies.
For an increasing dimension $D=2^6,2^7,\hdots,2^{14}$, we draw support sets
$T\in\Gamma$ of fixed sizes $M=4,8,16,32$ and random Fourier coefficients
(normal distribution for OMP and BP, modulus one for Thresholding).
We then test for the {\em smallest} number $N$ of continuously drawn samples,
such that at least $90$ percent ($180$ out of $200$) of the runs result in a
perfect recovery of the given Fourier coefficients.
Figure \ref{fig:ex4}(b-d) confirm the relation $\theta=C_{\epsilon}
\log_2(D)$ to reach a fixed success rate.
In contrast to Thresholding, the constant $C_{\epsilon}$ even decreases mildly
for OMP and BP when using a larger number $M$ of non-zero coefficients.

Both might indicate that the number of necessary samples is
$N=CM\log(D/(M\epsilon))$ rather than $N= C M \log(D/\epsilon)$ for OMP and BP
in an average case situation.
Note furthermore, that the Gaussian ensemble indeed allows for reconstruction
if this condition \eqref{Gauss_delta} is fulfilled.
\end{example}
\goodbreak

\begin{figure}[ht!]
  \centering
  \subfigure[Average computation time ($10$ runs) vs. dimension
    $D=2^3,2^4,\hdots,2^9$, $M=\lfloor\frac{1}{2}\sqrt{D}\rfloor$, and
    $N=M(\log_2(D)-2)$.
    Tested algorithms:
    OMP/QR/update (solid),
    BP/MOSEK (dashed),
    BP/CVX (dash-dot),
    BP/L1MAGIC/real/explicit (dotted), all with NFFT.
  ]{\includegraphics[width=0.45\textwidth]{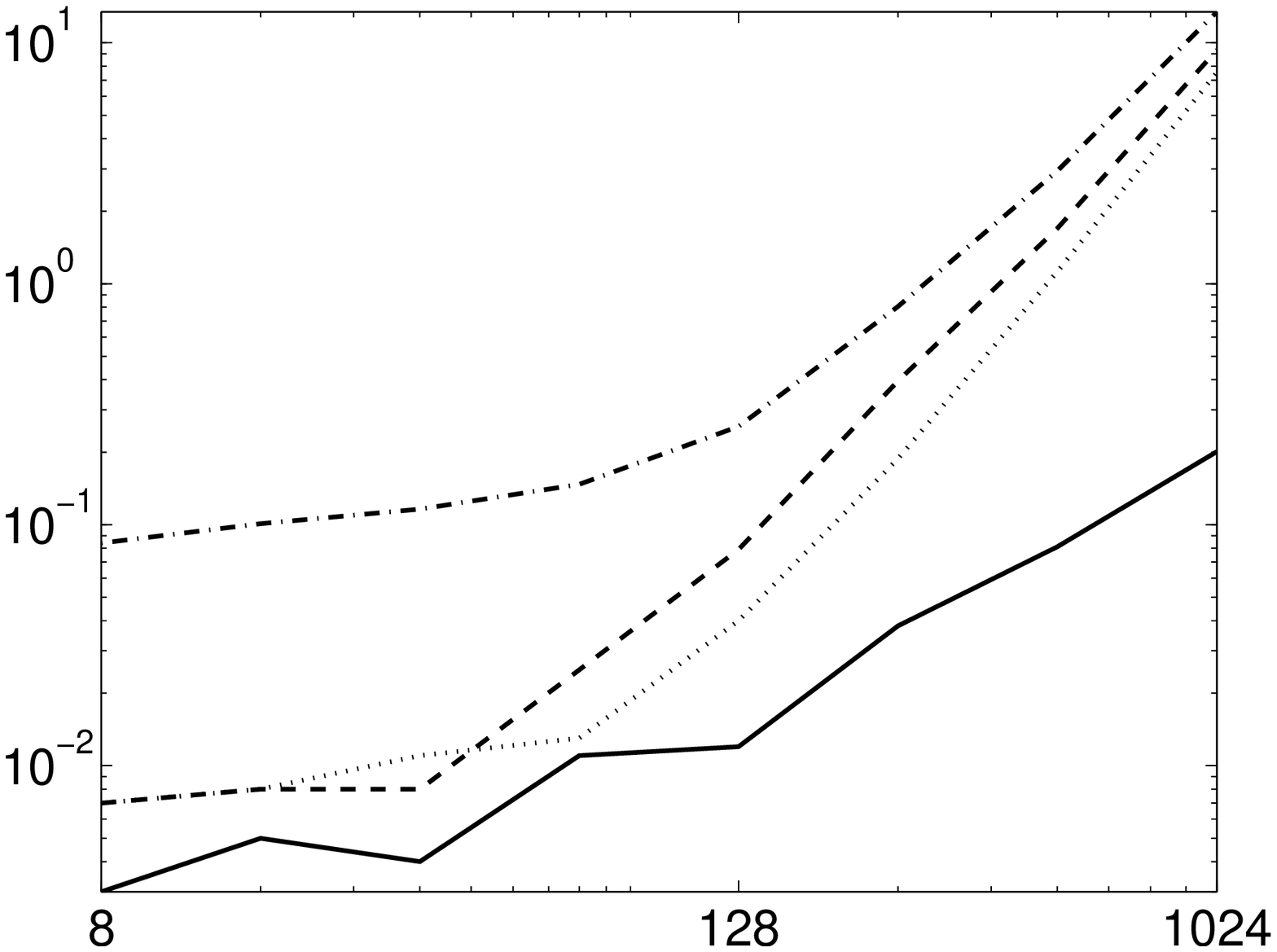}}\hfill
  \subfigure[Computation time vs. dimension $D=2^7,2^8,\hdots,2^{17}$,
    $M=\lfloor\frac{1}{8}\sqrt{D}\rfloor$, and $N=2M\log_2 D$.
    Tested algorithms:
    OMP/LSQR explicit (solid) and implicit (dashed);
    BP/L1MAGIC/real explicit (dash-dot) and implicit (dotted), all with FFT.
    For comparison: $10^{-7} D^{1.5} \log D$ (solid+diamond).
  ]{\includegraphics[width=0.45\textwidth]{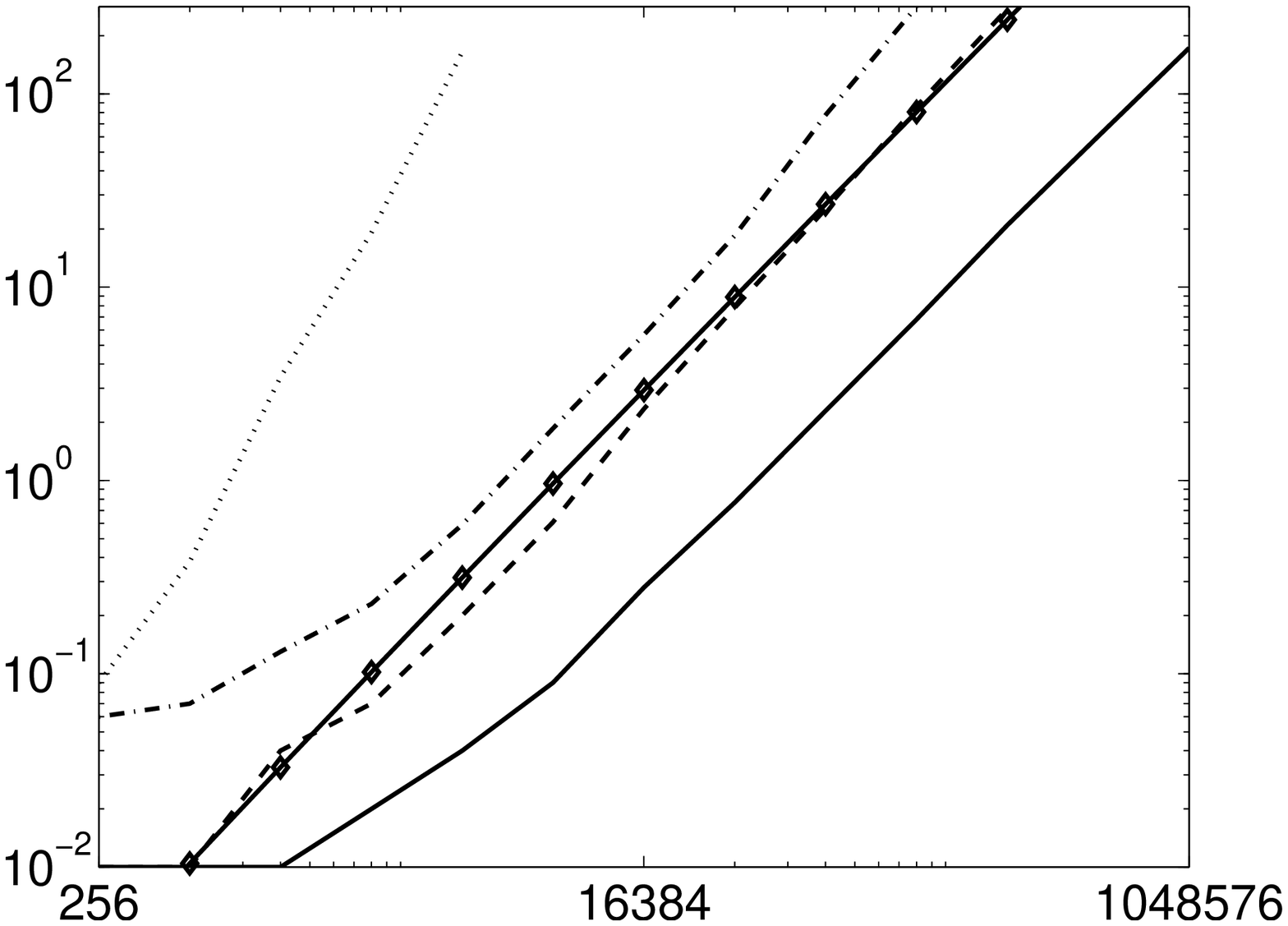}}\\
  \caption{Computation time in seconds with respect to an increasing dimension
    $D$ and dependent numbers $M$ and $N$ of non-zero coefficients and
    samples, respectively.\label{fig:ex7}}
\end{figure}

\begin{example}
% File    example7.m
This example considers the computation time needed by OMP and BP for an
increasing dimension $D$, a dependent number of non-zero coefficients
$M={\cal O}(\sqrt{D})$, and a number of samples $N={\cal O}(M\log D)$.
Constants are adjusted such that the used algorithms succeed in most cases in
the reconstruction task; all methods, except L1MAGIC, are tested with
complex coefficients.
For small dimensions $D$, we draw the samples continuously and test OMP with
the updated QR factorization and BP algorithms based on CVX, MOSEK, and
L1MAGIC.

Furthermore, we test with a somewhat larger dimension and discrete drawn
samples the algorithms:
OMP with LSQR and explicitly stored matrices $\F_{T_s X}$,
OMP with LSQR and FFT-based multiplications with $\F_{T_s X}$,
BP using L1MAGIC and an explicitly stored matrix $\F_X$ for $D\le 2^{12}$, and
BP using L1MAGIC and FFT-based multiplications with $\F_X$.
The FFT-based multiplications are denoted by implicit in Figure \ref{fig:ex7}
and need no additional memory, whereas storing the matrices explicitly needs
${\cal O}(D\log D)$ bytes for OMP and ${\cal O}(D^{1.5}\log D)$ bytes for BP.

Both OMP algorithms show a ${\cal O}(D^{1.5}\log D)$ time complexity, whereas
the scheme with explicit storage of $\F_{T_s X}$ is a constant multiple
faster.
The Basis Pursuit algorithms are considerably slower in all cases.
Moreover, the storage of the whole measurement matrix $\F_X$ results in large
memory requirements.
\end{example}

\begin{figure}[h]
  \centering
  \subfigure[Trig. polynomial and (noisy) samples.]   
  {\includegraphics[width=0.45\textwidth]{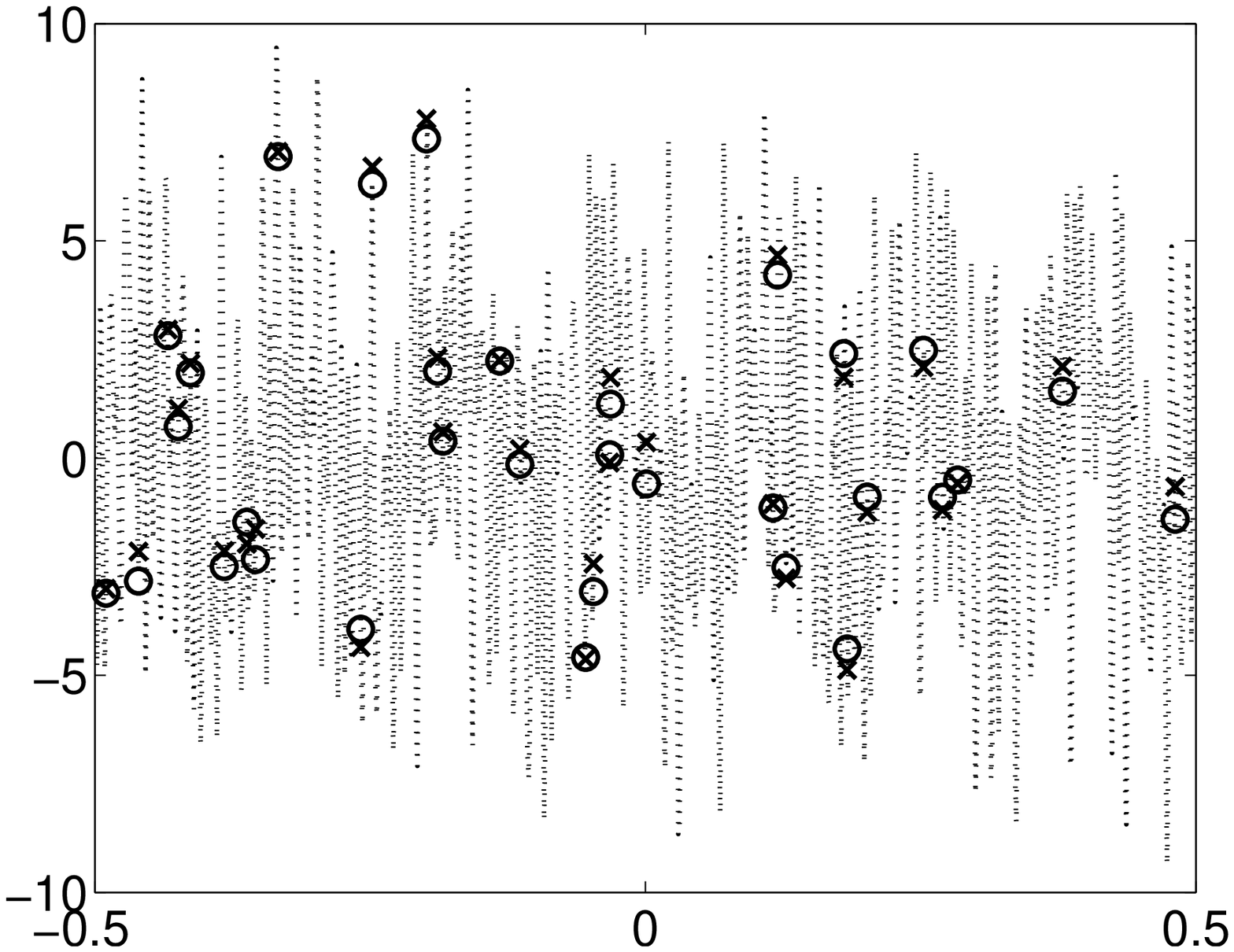}}\hfill
  \subfigure[True and recovered coefficients.]
  {\includegraphics[width=0.45\textwidth]{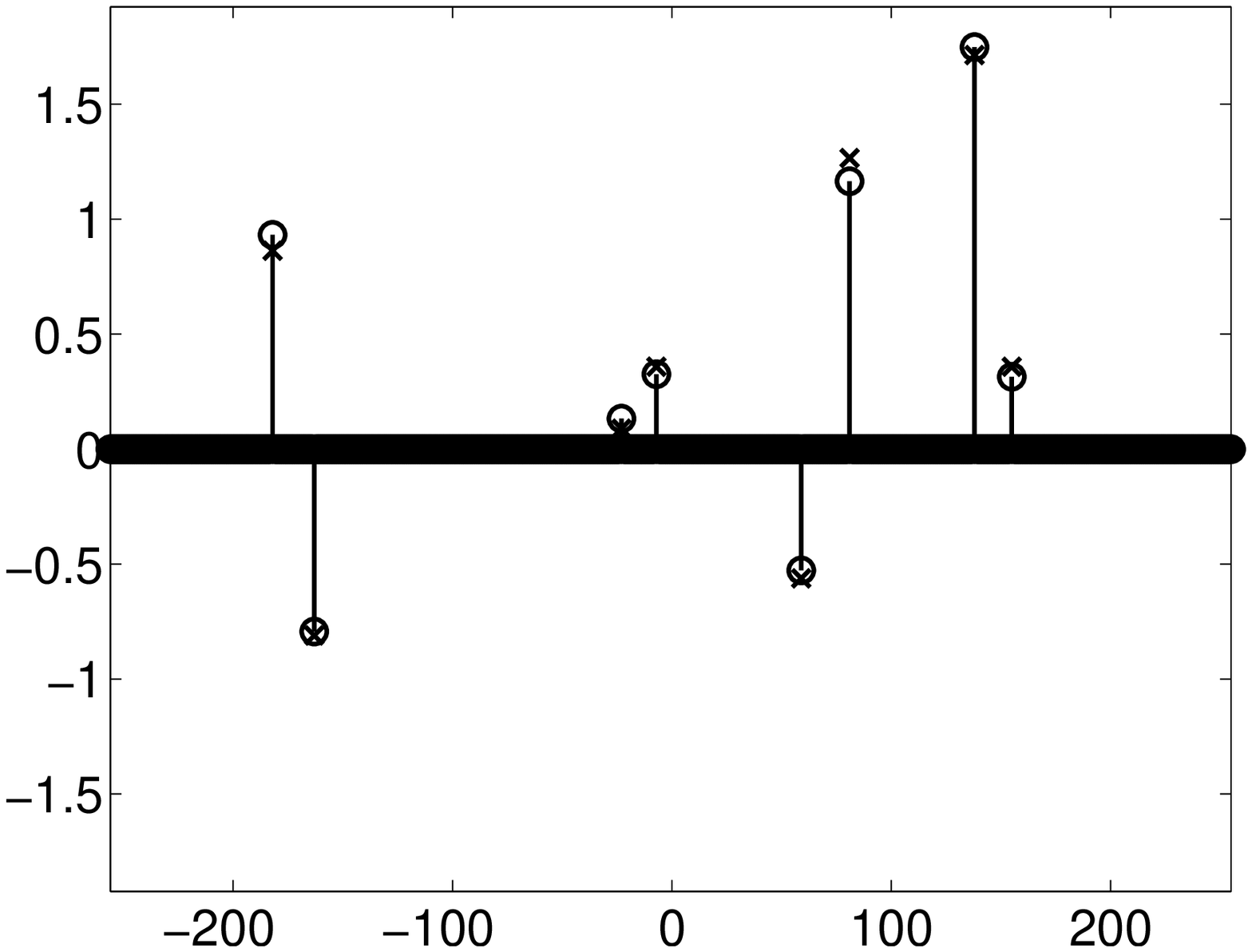}}\\
  \caption{Trigonometric polynomial (real part) as in Figure \ref{fig:1} and
    $30$ samples ($\circ$).
    The samples are disturbed by Gaussian distributed noise with 
    variance $0.2$ ($\times$) (resulting here in a PSNR of $15.6$dB). 
    Nevertheless OMP reconstructs the true support set of the coefficients, and
    the true coefficients ($\circ$) themselves with small error ($\times$).
    %, using 6 iterations instead of the necessary
    %three, 
    %the true coefficients ($\circ$) almost perfectly
    %($\times$).
    \label{fig:8}}
\end{figure}

\begin{example} 
% File: example8.m
Of course, stability under noise is important in practice.
For Basis Pursuit this has been established already by Cand{\`e}s, Romberg and
Tao in \cite{CRT2}.
For Orthogonal Matching Pursuit we performed first numerical experiments in
the following way.
We corrupted the sample values with a significant amount of normal
distributed noise. We observed that OMP usually 
finds the correct support set and makes only small errors 
on the coefficients, see Figure \ref{fig:8} and also \cite{Rau3}. 
Thus, it seems
that also OMP is stable under noise
 -- at least if the noise level is not too
high. (Actually, we observed that for moderately higher noise as in 
our example in Figure \ref{fig:8}, the reconstructed coefficient vector
is significantly different from the original). 
%Such issues are more thoroughly addressed in the follow-up
%paper \cite{Rau3}.
\end{example}

At last, note that both OMP and BP can be used to identify dominant
frequencies from few samples -- even if these frequencies are very high or if
two (high) neighboring frequencies are present.
Numerical tests showed for instance that the FFT based OMP can recover a
signal consisting of $10$ frequencies in $\{-2^{18},\hdots,2^{18}-1\}$ -- with
$k_1=2^{18}-1=262143$, $k_2=2^{18}-2=262142$ being two of them -- from $60$
random samples.

\section{Conclusions and future work}

Our theoretical and numerical results indicate that both Basis Pursuit and
greedy algorithms such as Orthogonal Matching Pursuit and Thresholding are
well-suited for the problem of recovering sparse trigonometric polynomials
from few random samples taken either on a grid (FFT) or on the cube (NFFT).

In practice, OMP obtains much better success rates than Thresholding.
Moreover, our numerical results indicate that recovery success rates at
reasonably small generalized oversampling rates $N/M$ are very similar on
generic signals for OMP and BP --although it is known that BP gives a uniform
recovery guarantee while greedy algorithms only provide a non-uniform
guarantee. 
As Theorem \ref{thm:omp} just analyzes the first step of OMP it is still
open to investigate theoretically the subsequent iterations. 
Due to stochastic dependency issues this does not seem to be straightforward,
see Remark \ref{rem:iterations}.

Both greedy methods are significantly faster than standard methods for Basis Pursuit.
In contrast, the series of papers \cite{GGIMS,GiMuSt,ZGSD,Zou} following
\cite{Mansour} considers algorithms that allow for asymptotic running time
proportional to $\log^{\beta}(D)$ for some $\beta$.
It would be very interesting to compare the two approaches numerically.
Moreover, step 4 of Algorithm \ref{algo:omp} might be replaced by so-called
isolation and group testing techniques as used in \cite{ZGSD,Zou}.

Typically, signals are not sparse in a strict sense but might still be
well-approximated by sparse ones.
Due to \cite{CRT2}, recovery by Basis Pursuit is still possible with small
errors provided that the restricted isometry constants are small.
As was noted in Section~\ref{Sec:BP}, our random Fourier matrices $\F_X$
satisfy this condition \cite{CT,RV06,Rau3}.

\bigskip
{\bf Acknowledgments:} 
The first author is grateful for partial support of this work by the German
Academic Exchange Service (DAAD) and the warm hospitality during his stay at
the Numerical Harmonic Analysis Group, University of Vienna.
The second author is supported by an Individual Marie Curie Fellowship funded
by the European Union under contract MEIF-CT-2006-022811. 
He also acknowledges funding
by the European Union's Human Potential Programme, under contract
HPRN--CT--2002--00285 (HASSIP).

We would like to thank Emmanuel Cand{\`e}s, Justin Romberg, Ingrid Daubechies,
Jared Tanner, Joel Tropp and Roman Vershynin for interesting discussions on the subject.

\begin{tabular}{ll}
Stefan Kunis & Holger Rauhut\\
Fakult\"at f\"ur Mathematik & NuHAG, Faculty of Mathematics\\
Technische Universit\"at Chemnitz & University of Vienna\\
Reichenhainer Str.~39 & Nordbergstr.~15\\
09107 Chemnitz, Germany & A-1090 Wien, Austria\\
kunis@mathematik.tu-chemnitz.de & holger.rauhut@univie.ac.at\\
\tt www.tu-chemnitz.de/$\sim$skunis &
\tt homepage.univie.ac.at/holger.rauhut
\end{tabular}
\end{document}